\documentstyle[11pt]{article}
\newcommand{\D}[2]{D^{#1}_{#2}}
\newcommand{\dubbelR}{{\sf I}\kern-.12em{\sf R}}
\newcommand{\dubbelN}{{\sf I}\kern-.12em{\sf N}}
\newcommand{\dubbelZ}{\mbox{\sf Z\hspace*{-5.4pt}Z}}
\newfont{\Bbb}{msbm10 scaled\magstep1}
\newcommand{\CC}{\hbox{\Bbb C}}

\def\unit{{\sf1 \! I}}
\def\half{{\frac{1}{2}}}
\newcommand{\ij}{i\kern -0.08em j}

\newcommand{\be}{\begin{equation}}
\newcommand{\ee}{\end{equation}}
\newcommand{\bea}{\begin{eqnarray}}
\newcommand{\eea}{\end{eqnarray}} 
\newcommand{\nn}{\nonumber}
\newcommand{\dubbelint}{\mathop{\int\int}\limits_{\kern-5.5pt {V-V_1}}}
\newcommand{\al}{\alpha}
\newcommand{\bt}{\beta}
\newcommand{\gm}{\gamma}

\newcommand{\dl}{\delta}
\newcommand{\Dl}{\Delta}
\newcommand{\ep}{\epsilon}

\newcommand{\th}{\theta}

\newcommand{\lm}{\lambda}

\newcommand{\sg}{\sigma}
\newcommand{\Sg}{\Sigma}

\newcommand{\om}{\omega}

\newcommand{\FSA}{\cal A}
\newcommand{\FSB}{\cal B}

\newcommand{\winv}{w^{-1}}
\newcommand{\xinv}{x^{-1}}
\newcommand{\yinv}{y^{-1}}
\newcommand{\zinv}{z^{-1}}
\newcommand{\ginv}{g^{-1}}
\newcommand{\hinv}{h^{-1}}
\newcommand{\ninv}{n^{-1}}
\newcommand{\minv}{m^{-1}}
\newcommand{\uv}{u^{-1}v}
\newcommand{\yx}{y(\xi)}
\newcommand{\yxinv}{y(\xi)^{-1}}
\newcommand{\wx}{w(\xi)}
\newcommand{\wxinv}{w(\xi)^{-1}}
\newcommand{\fox}{f \otimes x}

\newcommand{\rphi}{{}^{r}\!\phi}
\newcommand{\rD}{{}^{r}\!D} 
\newcommand{\DA}{{\cal D}({\cal A})}
\newcommand{\DG}{{\cal D}(G)}

\newcommand{\Valbe}{V_{\al}\otimes V_{\bt}}

\newcommand{\PAal}{\Pi^A_{\al}}
\newcommand{\PBbe}{\Pi^B_{\beta}}
\newcommand{\PAB}{\PAal \otimes \PBbe}

\newcommand{\PCga}{\Pi^C_{\gamma}}

\newcommand{\NCom}{\left(\hat{N}_C\right)_{\om}}
\newcommand{\Nabc}{N^{AB\gamma}_{\al\beta C}}

\newtheorem{thm}{Theorem}[section]
\newtheorem{lem}[thm]{Lemma}
\newtheorem{prop}[thm]{Proposition}
\newtheorem{cor}[thm]{Corollary}
\newtheorem{defn}[thm]{Definition}
\newtheorem{ass}[thm]{Assumption}

\newcommand{\Proof}{\noindent{\bf Proof}\quad}
\newcommand{\halmos}{\hbox{\vrule height0.28cm width0.01cm\vbox{\hrule height
 0.01cm width0.3cm \vskip0.26cm \hrule height 0.01cm width0.3cm}\vrule
 height0.28cm width 0.01cm}}
\newcommand{\eof}{\quad\halmos}
\newcommand{\Conj}{{\rm Conj}}
\newcommand{\bsl}{\backslash}

\newcommand{\lmAB}{\lambda_{AB}}
\begin{document}
\title{
\hfill{\normalsize UvA-WINS-Wisk.\ 97-14 }\\[-0.3 cm]
\hfill{\normalsize UvA-WINS-ITFA 97-44}\\[1.5 cm]
Tensor product representations of the quantum double of a compact group}
\author{T.H. Koornwinder\thanks{email: thk@wins.uva.nl}
      \\KdV Institute for Mathematics, University of Amsterdam, \\
        Plantage Muidergracht 24,
        1018 TV Amsterdam, The Netherlands.
 \and F.A. Bais\thanks{email: bais@phys.uva.nl}$\;$ and
N.M. Muller\thanks{email: nmuller@phys.uva.nl}
      \\Institute for Theoretical Physics, University of Amsterdam, \\
        Valckenierstraat 65,
        1018 XE Amsterdam, The Netherlands.}
\date{January 22, 1998.}
\maketitle
\begin{abstract}
We consider the quantum double $\DG$ of a compact group $G$,
following an earlier paper.
We use the explicit comultiplication on $\DG$ in order to build
tensor products of irreducible $*$-representations. Then we study
their behaviour under the action of the $R$-matrix, and their
decomposition into irreducible $*$-representations.
The example of ${\cal D}(SU(2))$ is treated in detail, with
explicit formulas for direct integral decomposition (`Clebsch--Gordan
series') and
Clebsch-Gordan coefficients. We point out possible physical applications.
\end{abstract}
\section{Introduction\label{s1}}
\setcounter{equation}{0}
Over the last decade quantum groups have become an important subject of
research both in mathematics and physics, see a.o.\ the monographs
\cite{CharPres}, \cite{Jant}, \cite{Kass} and \cite{Majbook}.
Of special importance are those quantum groups which are
quasi-triangular Hopf algebras, and thus have a universal $R$-element
satisfying the quantum Yang--Baxter equation. Via the QYBE there is a
connection with the braid group and thus with the theory of invariants
of links and 3-manifolds. In the physical context quantum groups 
play an important role in the theory of integrable lattice models,
conformal field theory (Wess--Zumino--Witten models for example) and
topological field theory (Chern--Simons theory).

Drinfel'd \cite{Drin} has introduced the notion of the quantum double
$\cal D(\cal A)$ of a Hopf algebra ${\cal A}$.
His definition (rigorous if $\cal A$
is finite dimensional, and formal otherwise) yields a
quasi-triangular Hopf algebra $\cal D(\cal A)$ containing ${\cal A}$
as a Hopf subalgebra. 
For $\cal A$ infinite dimensional, various rigorous definitions for the
quantum double or its dual have been proposed,
see in particular Majid \cite{Majbook} and Podles' and Woronowicz' \cite{PW}.

An important mathematical application of the Drinfel'd double is a
rather simple construction of the `ordinary' 
quasi-triangular quantum groups (i.e.\ $q$-deformations of universal
enveloping algebras of semisimple Lie algebras and of
algebras of functions on the corresponding groups), see for
example \cite{CharPres} and \cite{Majbook}. 

In physics the quantum double has shown up in various places:  
in integrable field theories \cite{BerLC}, in 
algebraic quantum field theory \cite{Mueger}, and in lattice quantum field 
theories. For a short summary of these applications, see
\cite{HausNill}. Another interesting application lies in orbifold
models of rational 
conformal field theory, where the physical sectors in the 
theory correspond to irreducible unitary representations of the
quantum double 
of a finite group. This has been constructed by Dijkgraaf, Pasquier
and Roche in \cite{DPR}. Directly related to the latter are the models of 
topological interactions between defects in spontaneously broken gauge 
theories in 2+1 dimensions. In
\cite{BWP} Bais, Van Driel and De Wild Propitius show that
the non-trivial fusion and braiding properties of the excited states in
broken gauge theories can be fully described by the representation theory
of the quantum double of a finite group. For a detailed treatment see 
\cite{Mark}.

Both from a mathematical and a physical point of view it is interesting to
consider the quantum double $\DG$ of the Hopf $*$-algebra
of functions on a (locally) compact group $G$, and to study its representation
theory.
For $G$ a finite group, $\DG$ can be realized as the linear space of
all complex-valued functions on $G\times G$. Its Hopf $*$-algebra
structure, which rigorously follows from Drinfel'd's definition, can be
given explicitly. In \cite{KM} and in the present paper we take the following
approach to $\DG$ for $G$ (locally) compact:\\
We realize
$\DG$ as a linear space in the form
$C_c(G\times G)$, the space of
complex valued, continuous functions of compact support on $G\times G$.
Then the Hopf $*$-algebra operations for $G$ finite can be  formally
carried over to operations on $C_c(G\times G)$ for $G$ non-finite
(formally because of the occurrence of Dirac delta's).
Finally it can be shown that these operations formally satisfy the axioms of
a Hopf $*$-algebra.

In \cite{KM},
we focussed on the $*$-algebra structure of
$\DG$, and we derived a classification of the irreducible
$*$-representations (unitary representations). In the present paper,
where we restrict ourselves to the case where $G$ is compact,
we address questions about `braiding' and `fusion' properties 
of tensor product representations of $\DG$,
for which the comultiplication and the $R$-matrix are explicitly
needed. We envisage physical applications in
nontrivial topological theories such as (2+1)-dimensional quantum
gravity, and higher dimensional models containing solitons \cite{BS}. 
In view of these and other applications we present our results
on representation theory not just abstractly, but quite explicitly.

\vspace*{.2cm}
The outline of the paper is as follows. 
In section \ref{s2} we specify the Hopf $*$-algebra structure of $\DG$. 
We then turn to the irreducible unitary representations in section
\ref{s10}, where we first recall a main result of \cite{KM}, concerning the
classification of these representations. We give a definition of their
characters, and compare the result to the case of finite $G$. 

An outstanding feature of quasi-triangular Hopf algebras is that their
non-cocommutativity is controlled by the $R$-element. Together with 
the explicit expression for the comultiplication 
this results in interesting properties of tensor products of irreducible
$*$-representations of $\DG$. In section \ref{s4} we define such a
tensor product representation, and specify the action of the quantum double. 
In section \ref{s5} we give the action of the universal $R$-matrix
on tensor product states (`braiding') on a formal level.
The rather non-trivial Clebsch--Gordan {\sl series} of irreducible 
$*$-representations (`fusion rules') are discussed in section \ref{s6}. 
They are calculated indirectly, namely, via direct projection of states, 
and the comparison of squared norms. This direct projection 
results in a very general method to construct the Clebsch--Gordan 
{\sl coefficients} of a quantum double in case orthogonal bases
can be given for the representation spaces. Finally, section \ref{s7}
treats the example of $G= SU(2)$ in detail.
\section{The Hopf algebra structure of $\DG$\label{s2}}
\setcounter{equation}{0}
Drinfel'd \cite{Drin} has given a definition of the {\sl quantum double} 
${\cal D}({\cal A})$ of a Hopf algebra ${\cal A}$. Write ${\cal A}^o$ for
the dual Hopf algebra to ${\cal A}$ with opposite comultiplication. Then
$\DA$ is a quasi--triangular Hopf algebra, it is equal to ${\cal A}
\otimes {\cal A}^o$ as a linear space, and it contains ${\cal A}\otimes 1$
and $1\otimes{\cal A}^o$ as Hopf subalgebras. If ${\cal A}$ is moreover
a Hopf $*$-algebra then $\DA$ naturally becomes a Hopf $*$-algebra. This
definition of the quantum double is only rigorous if ${\cal A}$ is finite
dimensional. 

If $G$ is a compact group and $C(G)$ the Hopf $*$-algebra of continuous
complex values functions on $G$, then instead of ${\cal D}(C(G))$ we will 
write $\DG$ for the quantum double of $C(G)$. For $G$ a finite group we have
\be
\DG \simeq C(G)\otimes \CC[G] \simeq C(G\times G)
\ee
as linear spaces. Also in the case of a finite group it is possible to
write down the formulas for the Hopf $*$-algebra operations and the
universal $R$-element of $\DG$, both in the formulation with $\DG =
C(G)\otimes \CC[G]$ (see \cite{DPR}) and with $\DG = C(G\times G)$. 
In the last picture
the formulas may typically involve a summation over the group or
a Kronecker delta on $G$. They suggest analogous formulas for $G$
arbitrarily compact, by simply replacing the summation over $G$ by
integration w.r.t.\ the normalised Haar measure on $G$, and replacing 
the Kronecker delta by the Dirac delta. This way we obtain the following
definitions, where $F, F_1, F_2 \in C(G\times G)$ and $x,y,x_1, y_1, x_2,
y_2 \in G:$\\
Multiplication:
\be
(F_1\bullet F_2)(x,y):=\int_G F_1(x,z)\,F_2(z^{-1}xz,z^{-1}y)\,dz.\\
\label{eq:multiplication}
\ee
$*$-operation:
\be
F^{*}(x,y) = \overline{F(\yinv xy,\yinv)}.
\label{eq:star}
\ee
Unit element
\be
1(x,y) = \delta_e(y),
\label{eq:unit}
\ee
Comultiplication:
\be
(\Dl F) (x_1, y_1; x_2, y_2) = F(x_1 x_2, y_1)\, \dl_e(y_1^{-1}y_2)
\label{eq:comult}
\ee
Counit:
\be
\ep (F) = \int_G F(e,y)\;dy
\ee 
Antipode 
\be
(S(F))(x,y) = F(\yinv\xinv y, \yinv)
\label{eq:antipode}
\ee
Universal $R$-element:
\be
R(x_1, y_1;x_2, y_2) = \dl_e(x_1\yinv_2)\,\dl_e(y_1).
\label{eq:Relement}
\ee

Note that due to the occurring Dirac delta's the unit element in fact 
does not lie inside $\DG$. Similarly, the comultiplication doesn't map
into $\DG\otimes \DG$ (not even into the topological completion
$Cl (C(G\times G)\otimes C(G\times G)) \simeq C(G\times G\times G\times G)$),
and furthermore the $R$-element doesn't lie inside $\DG\otimes \DG$. In 
practice this does not pose a serious problem as 
we will always formally integrate over these Dirac delta's,
nevertheless we still have to be careful in dealing with the resulting
expressions, because it can happen that the Dirac delta is partially
fulfilled, giving rise to infinities. 
 
With the above operations $C(G\times G)$ formally becomes a quasi-triangular
Hopf $*$-algebra called $\DG$. For the case of a finite group $G$ this holds
rigorously, which is clear just by the quantum double construction. However,
for the case of general compact $G$ we have to verify that
Eqs.(\ref{eq:multiplication})--(\ref{eq:Relement}) do indeed satisfy all
axioms of a quasi-triangular Hopf $*$-algebra. 

In \cite{KM} it was observed that $C(G\times G)$ with
Eqs.(\ref{eq:multiplication}) and (\ref{eq:star}) is a $*$-algebra, and
furthermore the irreducible unitary representations of 
this $*$-algebra were studied and classified. In the present paper we will
consider tensor products and braiding properties of these irreducible
$*$-representations (from now on mostly referred to as `irreps') by using
the comultiplication and the $R$-element.
\section{Irreducible representations \label{s10}}
\setcounter{equation}{0}
We recapitulate the contents of Corollary 3.10, one of the main results 
of \cite{KM}. 
Throughout, when we speak of a compact group (or space), we 
tacitly assume that it is a separable compact Hausdorff group (or space).
\begin{defn}\label{def:rep}
Let $G$ be a  compact group, and $\Conj(G)$ the 
collection of conjugacy classes of $G$ (so the elements of $\Conj(G)$ are
the sets of the form $\{xg\xinv\}_{x\in G}$ with $g\in G$).  For each
$A\in \Conj(G)$ choose some representative $g_A\in A$, and let $N_A$ be
the centralizer of $g_A$ in $G$. For each $\al \in
\widehat{N_A}$ (the set of equivalence classes of irreducible
unitary representations of $N_A$) choose
a representative, also denoted by $\al$, which is an irreducible unitary
representation of $N_A$ on some finite dimensional Hilbert space $V_{\al}$. 
Also, let $dz$ be the normalised Haar measure on $G$. 
For measurable functions $\phi:G\to V_{\al}$ such that for all $h\in N_A$
it holds that
\be
\phi(gh) = \al(\hinv)\phi(g)\quad\mbox{ for almost all }\,g\in G,
\label{eq:covariance}
\ee
we put
\be
\|\phi\|^2 := \int_{G} \|\phi(z)\|^2_{V_{\al}}\,dz.
\label{eq:repfunctions}
\ee
Now $L^2_{\al}(G,V_{\al})$, which is the linear space of all such $\phi$
for which $\|\phi\| <\infty$ divided out by the functions with norm zero,
is a Hilbert space. 
\end{defn}
The elements of $L^2_{\al}(G,V_{\al})$ can also be considered as 
$L^2$-sections of a homogeneous vector bundle over $G/N_A$. The space
$L^2_{\al}(G,V_{\al})$ is familiar as the representation space of the
representation of $G$ which is induced by the representation $\al$ of
$N_A$.
\begin{thm}\label{thm:rep}
For $A\in \Conj(G)$ and $\al\in\widehat{N_A}$
we have mutually inequivalent
irreducible $*$-representations $\PAal$ of $\DG = C(G\times G)$ on $L^2_{\al}
(G,V_{\al})$ given by
\be
\left( \PAal (F) \phi\right)(x) := \int_G\, F(xg_A\xinv, z) \phi(\zinv x)\,dz,
\quad F\in \DG
\label{eq:representation}
\ee
These representations are moreover $\|.\|_1$-bounded (see for this notion
formula (33) in \cite{KM}).
All irreducible $\|.\|_1$-bounded $*$-representations of $\DG$ are equivalent
to some $\PAal$.
\end{thm}

In fact, a much more general theorem holds (see Theorem 3.9 in \cite{KM}),
namely for the representation theory of so-called transformation group 
algebras $C(X\times G)$, where the compact group $G$ acts continuously
on the compact set $X$, instead of the conjugation action of $G$ on $G$. 
We may even assume $G$ and $X$ to be locally compact, under the extra 
condition of {\sl countable separability} of the $G$-action. Then we have 
to consider $C_c(X\times G)$ and use a quasi-invariant measure on $G/N_A$.

Also, the rest of the Hopf algebra structure of $\DG$, in particular
the comultiplication, will survive for the case of noncompact $G$ as
long as $G$ acts on itself by conjugation. It would be interesting
to extend the results of this paper to this case of (special) noncompact
$G$.

An interesting issue in representation theory is the character  
of an irrep. For the case of a finite group $G$ such characters have been
derived in \cite{DPR}. For our case, where irreps are generally infinite
dimensional, the operator $\PAal(F)$ will not be trace class for all
$F\in \DG$, so we restrict ourselves to the case of a Lie group
$G$ and $C^{\infty}$-functions on $G\times G$. In this paper we will
only state the formula for the characters of irreps
of the quantum double. The proof for it, the orthogonality of the
characters, and the related
subject of harmonic analysis, will be given in a forthcoming paper.
\begin{thm}\label{thm:char}
Let $\chi_{\al}$ denote the character of the irreducible
$*$-representation  $\al$ of $N_A$.
For an irreducible $*$-representation $\PAal$ of the quantum double $\DG$
the character is given by
\be
\chi^A_{\al}(F) = \int_G\,\int_{N_A}\, F(z g_A \zinv, zn \zinv)\,
\chi_{\al}(n)\,dn\,dz, \qquad F\in C^{\infty}(G\times G).
\label{eq:char}
\ee
\end{thm}

Let us check the connection with the case of a finite group $G$. 
As discussed in \cite{KM} for a finite group $G$ there is a
linear bijection $\DG = C(G)\otimes\CC[G] \Longleftrightarrow C(G\times G)$:
\bea
\fox &\mapsto& \left((y,z) \mapsto f(y) \dl_x(z)\right) \nn\\
\sum_{z\in G} F(.\, ,z)\otimes z & \leftarrow& F
\eea
Taking $f = \delta_g$ as function on (finite) $G$, we obtain
\be
\chi^A_{\al} (\delta_g\otimes x) = \int_G\int_{N_A} \delta_g(zg_A\zinv)
\delta_x(zn\zinv)\chi_{\al}(n)\,dn\,dz,
\ee
which indeed coincides with the definition of the character in \cite{DPR}.
\section{Tensor products \label{s4}}
\setcounter{equation}{0}
In section \ref{s10} we have recapitulated the classification of the 
irreducible
$*$-representations of the quantum double $\DG$. With the coalgebra structure 
of $\DG$ that we have derived in  section \ref{s2} we
can now consider tensor products of such  representations.

Let $\PAal$ and $\PBbe$ be irreducible $*$-representations of $\DG$.
For the representation space of the tensor product representation we
take the Hilbert space of vector-valued functions on $G\times G$ as
follows: for measurable functions $\Phi\,:\,G\times G \to V_{\al}
\otimes V_{\bt}$ such that for all $h_1\in N_A, h_2\in N_B$ it holds that
\be 
\Phi(x h_1, y h_2) = \al (\hinv_1)\otimes \bt (\hinv_2) \Phi (x, y),
\qquad\mbox{for almost all}\; (x,y) \in G\times G
\label{eq:tenscov}
\ee
we put
\be
\|\Phi\|^2 := \int_G\int_G\| \Phi (x,y)\|^2_{V_{\al}
\otimes V_{\bt}} dx\, dy.
\label{eq:tensnorm}
\ee
Now the space $L^2_{\al\bt}(G\times G, V_{\al}\otimes V_{\bt})$ is
defined
as the linear space of all such $\Phi$ for which $\|\Phi\|<\infty$,
divided out by the functions of norm zero. 
Note that this space is the completion of
the algebraic tensor product of $L^2_\al(G,V_\al)$ and $L^2_\bt(G,V_\bt)$.
By Eq.(\ref{eq:representation}) the tensor product representation 
$\PAal \otimes \PBbe$ becomes formally:
\bea
\left(\left(\PAal \otimes \PBbe\right)(F)\,\Phi\right)(x, y) &:=& 
\left( \left(\PAal \otimes \PBbe\right)(\Delta F)\,\Phi\right)(x,y)\nn\\
&=&\int_G \int_G\,\Delta F(x g_A \xinv, z_1; y g_B \yinv, z_2)\, 
\Phi(\zinv_1 x, \zinv_2 y)\,dz_1\,dz_2. \nn
\eea
Then it follows by substitution of Eq.(\ref{eq:comult}) and by
formally integrating the Dirac delta function that
\bea
\left( \left(\PAal \otimes \PBbe\right)(F)\,\Phi\right)(x, y)
&=&\int_G F(x g_A \xinv y g_B \yinv, z)\,\Phi(\zinv x, \zinv y)\,dz.
\label{eq:tensproduct}
\eea
It is easy to see that this
is indeed a representation of $\DG$: there is the covariance property, 
as given in Eq.(\ref{eq:tenscov}), and the homomorphism 
property can be readily checked. The functions of the form
\be
\Phi(x, y) = \phi^A_{\al}(x)\otimes \phi^B_{\bt}(y) \in V_{\al}
\otimes V_{\bt}
\ee
(with $\phi^A_{\al}$ and $\phi^B_{\bt}$  basis functions of 
the representation spaces for $\PAal$ and $\PBbe$ respectively) span a 
dense subspace of $L^2_{\al\bt}(G\times G,V_{\al}\otimes V_{\bt})$. 
The positive--definite inner product then reads
\be
\langle \Phi_1, \Phi_2\rangle := \int_{G}\int_{G}
\langle{\phi_1}^A_{\al}(x),{\phi_2}^A_{\al}(x)\rangle_{V_{\al}}\,
\langle{\phi_1}^B_{\bt}(y),{\phi_2}^B_{\bt}(y)\rangle_{V_{\bt}}\,
dx\,dy.
\ee
This tensor product representation now enables us to further analyse 
two important operations which are characteristic for quasi-triangular 
Hopf algebras, namely `braiding' and `fusion'. They will turn up in several
applications of these algebras \cite{BM}.
\section{Braiding of two representations\label{s5}}
\setcounter{equation}{0}
Let us investigate  
the action of the universal $R$-element in the aforementioned tensor product
representation. A simple formal calculation with use of
Eqs.(\ref{eq:Relement}) and (\ref{eq:representation}) yields
\bea
\left( \left(\PAal \otimes \PBbe\right)(R)\,\Phi\right)(x, y) &=& 
 \int_G \int_G \, \delta_e(x g_A\xinv\zinv)\,\delta_e(w)\,\Phi(w^{-1}x,
 \zinv y)\, dw\,dz \nn \\
&=& \Phi(x, x g_A^{-1}\xinv y).
\eea
The braid operator $\cal{R}$ is an intertwining mapping between
$\PAal \otimes \PBbe$ on $V_{\al}\otimes V_{\bt}$ and $\PBbe\otimes\PAal$
on $V_{\bt}\otimes V_{\al}$ given by
\be
{\cal{R}}^{AB}_{\al\bt}\,\Phi := \left(\sigma_L \circ \left(\PAal 
\otimes \PBbe \right)(R)\right) \Phi
\ee
where 
\be
(\sigma_L  \Phi)(x,y) := \sigma \left(\Phi(y,x)\right),\quad
\sigma (v\otimes w) := w\otimes v,\;\;v\in V_{\al}, w\in V_{\bt},
\ee
so it interchanges the representations $\PAal$ and $\PBbe$. Hence
\be
\left({\cal{R}}^{AB}_{\al\bt}\,\Phi\right)(x,y) = \left(\sigma_L \left( 
\left(\PAal \otimes \PBbe\right)
(R) \Phi(x,y)\right)\right) = \sigma \left(\Phi(y, y g_A^{-1}\yinv x)\right).
\label{eq:rab}
\ee
To make sure that Eq.(\ref{eq:Relement}), being  derived from
a formally defined $R$-element Eq.(\ref{eq:Relement}), yields the
desired intertwining property for ${\cal{R}}^{AB}_{\al\bt}$, one can derive
this property directly from
Eqs.(\ref{eq:rab}) and (\ref{eq:tensproduct}). Then we must show that
\be
\left({\cal{R}}^{AB}_{\al\bt}\left(\PAB\right) (F)\, \Phi\right)(x,y)
= \left(\left(\PBbe\otimes\PAal\right)
(F)\left({\cal{R}}^{AB}_{\al\bt}\,
\Phi\right)\right)(x,y).
\label{eq:Rcheck}
\ee
The right hand side of this equation gives
\bea
\lefteqn{\int_G F(x g_B \xinv y g_A\yinv, z)\left({\cal{R}}^{AB}_{\al\bt}\,
\Phi\right)(\zinv x, \zinv y)\,dz =\qquad\qquad}\nn\\
& & \int_G F(xg_B\xinv yg_A\yinv, z)
 \,\sg\left( \Phi(\zinv y, \zinv y\ginv_A\yinv x)\right)\, dz.
\eea
which is obviously equal to the left hand side of Eq.(\ref{eq:Rcheck}),
using Eq.(\ref{eq:rab}) and Eq.(\ref{eq:tensproduct}). 
\section{Tensor product decomposition   \label{s6}}
\setcounter{equation}{0}
Another general question is the decomposition of the tensor
product of two irreducible representations into irreducible 
representations:
\be
\PAal \otimes \PBbe \simeq \bigoplus_{C,\gm} \Nabc \,\PCga,
\label{eq:fusion}
\ee 
where we suppose that such a tensor product is always reducible. 
For finite $G$ tensor products of irreps of $\DG$ indeed decompose into
a direct sum over single irreps. For compact $G$ the direct sum over 
the conjugacy class label $C$ has to be replaced by a direct integral, 
\be
\PAal \otimes \PBbe  \simeq \bigoplus_{\gm} \int^{\oplus}\Nabc
\,\PCga\, d\mu(C) 
\label{eq:fusionintegral}
\ee
where $\mu$ denotes an {\sl equivalence class of measures} on the set of 
conjugacy classes, but the multiplicities must be the same for different
measures in the same class, see for instance the last Conclusion in
\cite{Arv} for generalities about direct integrals.   
Recall that two (Borel) measures $\mu$ and $\nu$ are equivalent iff they 
have the same sets of measure zero \cite{Arv}. 
By the Radon--Nikodym theorem, $\mu$ and $\nu$ are equivalent iff $\mu
= f_1 \nu,\; \nu = f_2 \mu$ for certain measurable functions $f_1, f_2
\geq 0$. If one considers specific states and/or
their norms (so elements of specific Hilbert spaces), it is required to 
make a specific choice for the measure. But if one only compares equivalence
classes of irreps, like we do in the 
Clebsch--Gordan series in Eq.(\ref{eq:fusionintegral}), the exact 
measure on $\Conj(G)$ is not of importance, only its equivalence class. 

Our aim is to determine the measure $\mu$ (up to equivalence) and the
multiplicities $\Nabc$ of this Clebsch--Gordan
series for $\DG$. In physics these $\Nabc$ are often referred to as 
`fusion rules', as for example in \cite{DPR} for the case of $G$ a
finite group.  

In ordinary group theory the multiplicities can be determined using the
characters of representations. Recall that for a continuous group $H$ with
irreducible representations $\pi^a, \pi^b, \pi^c,...$ and
characters $\chi^a, \chi^b, \chi^c$  the number of times 
that $\pi^c$ occurs in the $\pi^a\otimes \pi^b$ is given by
\be
n^c_{ab} = \int_{h\in H} \overline{\chi^c(h)} \chi^a(h) \chi^b(h)\,dh.
\ee
Thus a direct computation of the multiplicities requires an integration
over the group. For the quantum double this approach is not very
attractive, and we have to take an alternative route. Furthermore, the
direct decomposition of the character of a tensor product of irreps into
a direct sum/integral over characters of single irreps is problematic,
since the tensor product character is not trace class,
while the single characters are. 

The rigorous approach we will take is to look at the decomposition in 
more detail, in the sense that we consider the projection of a {\sl
state} in $\PAB$ onto 
{\sl states} in the irreducible components $\PCga$. Subsequently we 
compare the squared norm of the tensor product state with the direct 
sum/integral of squared norms of the projected states. This will lead to 
an implicit equation for the multiplicities $\Nabc$. 
The projection involves the construction of intertwining operators from the 
tensor product Hilbert space to Hilbert spaces of irreducible representations.
This construction is described in the next subsection, and the intertwining
operators are given 
in Theorem \ref{thm:rho}. If orthonormal bases are given for the 
Hilbert spaces of irreducible representations this means we can
derive the Clebsch--Gordan {\sl coefficients} for the quantum double.
In section \ref{s7} we will work this out explicitly for the case $G=SU(2)$.

Since the proof of Theorem \ref{thm:rho} is quite lengthy,
in the following paragraph we first give a brief outline of the
procedure we will follow.
\vspace*{.2cm}

To prove isometry between the Hilbert space of a tensor product
representation and a direct sum of Hilbert spaces of irreducible
representations we must construct an intertwining mapping $\rho$ from the
first space, whose elements are functions of two variables with a certain
covariance property, to the second space (=direct sum of spaces), whose 
elements are functions of one variable with a similar covariance property.
{}From Eq.(\ref{eq:tensproduct})
one can see that the conjugacy class label $C$ of the representation to which 
$\Phi$ must be mapped depends on the `relative difference' $\xi$ between the 
entries $(y_1,y_2)$ of $\Phi$. This $\xi$ is the variable that remains if
$(y_1,y_2)$ and $(z y_1 n_1, z y_2 n_2)$ are identified for all $z\in G$
and for all $n_1\in N_A$ and $n_2\in N_B$. So $\xi$ is an
element of the double coset $ G_{AB}=N_A\bsl G/N_B$ we have introduced
before. Eq.(\ref{eq:xyz}) in
Proposition \ref{prop:maps} shows how $C$ depends on $\xi$. 

In Proposition \ref{propmaps} we give a map ${\cal F}_1$ which constructs
a function $\phi$ on $G$ out of a function $\Phi$ on $G\times G$. The
action of $\DG$ on $\phi$ depends on the possible `relative differences'
$\xi$ between the entries of $\Phi$, which is why we say that
$\phi$ also depends on $\xi$. Therefore we introduce the function spaces 
of Eqs.(\ref{eq:funalbe}) and (\ref{eq:funom}).

Lemma \ref{lemone} shows that the squared norm of $\Phi$ equals the 
direct integral over $\xi$ of the squared norm of $\phi$, and thus that
the map ${\cal F}_1$ is an isometry of Hilbert spaces. One can also
think of $\xi$ as a {\sl label} on $\phi$, which distinguishes its behaviour
under the action of $\DG$, which is in fact shown by Lemma \ref{lemtwo}.
These two lemmas together provide the map $\PAB \to \int^{\oplus} 
\Pi^{C(\xi)}_{\om}\,d\mu(\xi)$ from the tensor product representation to
a direct integral over `single' (not yet irreducible) representations.

Subsequently we must decompose these representations $\Pi^{C(\xi)}_{\om}$
into irreducible representation $\PCga$. Comparing the covariance 
properties before and after $\rho$ we find the restriction on the
set $\gm$ may be chosen from, which is given in Eq.(\ref{eq:setNCom}).

Eq.(\ref{eq:G1p}) gives the isometry of a Hilbert space from the direct
integral of Hilbert spaces we constructed before (via the map ${\cal F}_1$)
into the direct sum of Hilbert spaces of irreducible representations $\PCga$.

The combination of these two steps in the tensor product decomposition
is summarised in Theorem \ref{thm:rho}. Finally we compare the squared
norms before and after the mapping $\rho$, and arrive at Eq.(\ref{117}),
which gives us an implicit formula for the multiplicities. 
The degeneracy of the irreducible representation
$\PCga$ depends on two things: firstly, the possible non-injectivity of the 
map $\xi \mapsto C$, which is taken into account by the integration over 
$N_A\bsl G/N_B$ with measure $dp^C(\xi)$. 
And secondly by the dimension $d_{\gm}$ of $V_{\gm}$. 
We now turn to the explicit proof.
\vspace*{.2cm}

To start with, fix the conjugacy classes $A$ and $B$, and also 
the irreducible unitary representations $\al\in\hat{N}_A$ and
$\bt\in\hat{N}_B$ with representation spaces $V_{\al}$ and $V_{\bt}$ 
of finite dimensions $d_{\al} = \mbox{dim}\,V_{\al}$ and  
$d_{\bt} = \mbox{dim}\,V_{\bt}$ respectively.
The set $\Conj(G)$ of conjugacy classes of $G$ forms a partitioning of $G$.
Therefore it can be equipped with the quotient topology, which is again
compact Hausdorff and separable. In Definition \ref{def:rep} we had already
chosen some representative $g_A \in A$ for each $A\in \Conj(G)$. We will
need the following assumption about this choice:
\begin{ass}
The representatives $g_A\in A$ can be chosen such that the map $A\mapsto
g_A\,:\,\Conj(G) \to G$ is continuous.
\label{ass:conjcont}
\end{ass}
In fact, we will make this particular choice. The assumption means that
the map from $G$ to $G$, which assigns to each $g\in G$ the representative
in its conjugacy class, is continuous. For $G$ a compact connected Lie
group we can make a choice of representatives $g_A$ in agreement with
Assumption \ref{ass:conjcont} as follows. Let $T$ be a maximal torus
in $G$, let $T_r$ be the set of regular elements of $T$ (i.e.\ those
elements for which the centraliser equals $T$), let
$K$ be a connected component of $T_r$, and let $\overline{K}$ be the
closure of $K$ in $T$. Take $g_A$ to be the unique element in the
intersection of the conjugacy class $A$ with $\overline{K}$. See for
instance reference \cite{BrockDieck}. 

Define
\be
G_{AB} := N_A\backslash G/N_B
\label{eq:doubcoset}
\ee
to be the collection of double cosets of the form $N_A y N_B,\; y\in G$. 
Then $G_{AB}$ also forms a partitioning of $G$ which can be equipped with
the quotient topology from the action of $G$ (compact Hausdorff and 
separable). Now also choose for each
$\xi \in G_{AB}$ some representative $y(\xi) \in \xi$. 
We will need the following assumption for this choice of representative:
\begin{ass}\label{ass2}
The representatives $y(\xi)\in \xi$ can be (and will be) chosen such that 
the map $\xi \mapsto y(\xi)\,:\, G_{AB} \to G$ is continuous.
\end{ass}
In other words, the map from $G$ to $G$ which assigns to each $g\in G$
the representative in the double coset $N_A g N_B$ is continuous.
For $SU(2)$ a choice of representatives $y(\xi)$ in agreement with
Assumption \ref{ass2} will be given in section \ref{ssseries}.

For $\xi \in G_{AB}$ define the conjugacy class $C(\xi)\in \Conj(G)$ by
\be
g_A y(\xi) g_B y(\xi)^{-1} \in C(\xi).
\label{eq:class}
\ee
Then the map
\be
\lmAB\,:\, \xi \mapsto C(\xi)\, : \, G_{AB} \to \Conj(G)
\label{eq:cosettoclass}
\ee
is continuous. Note that the image of $\lm_{AB}$ depends on the values of
$A$ and $B$,
but that $G_{AB}$ only depends on $N_A$ and $N_B$, so not on the precise
values of the conjugacy class labels. 

\begin{prop}\label{prop:maps}
(a)  We can choose a Borel map $\xi \mapsto \wx\,:\,G_{AB}\to G$ such that
\be
g_A \yx g_B \yxinv = \wx g_{C(\xi)} \wxinv
\label{eq:xyz}
\ee
(b)  We can choose a Borel map $x\mapsto (n_1(x), n_2(x))\,:\,G \to N_A
\times N_B$ such that
\be
x = n_1(x) y(N_A x N_B) n_2(x)^{-1}
\label{eq:borelns}
\ee
\end{prop}
\Proof
(a) The map 
\be
(w,C) \mapsto wg_C \winv \, :\, G\times \Conj(G) \to G \nonumber
\ee
is continuous (by Assumption \ref{ass:conjcont}) and surjective. By 
Corollary \ref{110} there exists a 
Borel map $x\mapsto (w_x, C_x) \,:\, G\to G\times \Conj(G)$ such that
$x = w_x g_{C_x} \winv_x $. Now take $x= g_A \yx g_B \yxinv$ for $\xi 
\in G_{AB}$, then from Eq.(\ref{eq:class}) it follows that $C_x = C(\xi)$.
The map $\xi \mapsto x$ is continuous by Assumption \ref{ass2}, the map 
$x \mapsto w_x$ is Borel.
Put $\wx := w_x$, then $\xi\mapsto \wx$ is Borel, and Eq.(\ref{eq:xyz})
is satisfied.\\
(b) The map 
\be
(n_1, n_2, \xi) \mapsto n_1 \yx \ninv_2\,:\,N_A\times N_B\times G_{AB}
\to G \nonumber
\ee
is continuous (by Assumption \ref{ass2}) and surjective.  By 
Corollary \ref{110} there exists a 
Borel map $x\mapsto (n_1(x), n_2(x), \xi(x))\,:\,G\to N_A\times N_B\times
G_{AB}$ such that $x=n_1(x) y(\xi(x)) n_2(x)^{-1}$. Then $\xi(x) = N_A
x N_B$, and thus Eq.(\ref{eq:borelns}) is satisfied.$\qquad$\eof 
  
\vspace*{.2cm}
Let $Z$ be the center of $G$, then $Z\subset N_A$ and $Z\subset N_B$.
By Schur's lemma $\al(z)$ and $\bt(z)$ will be a scalar for $z\in Z$.
Define the character $\om$ of $Z$ by
\be
\al(z)\otimes \bt(z) =: \om(z) \mbox{id}_{V_{\al}\otimes
V_{\bt}},\quad z\in Z. 
\label{eq:omega}
\ee 
With this character we now define the linear spaces
\bea
\mbox{Fun}_{\al,\bt}(G\times G,V_{\al}\otimes V_{\bt}) &:=& \{\Phi : G\times G
\to \Valbe \;|\;\Phi(u\ninv_1, v\ninv_2) = \al(n_1)\otimes\bt(n_2)\, 
\Phi(u,v) \nn \\
& &\forall n_1\in N_A, n_2\in N_B,\; u,v\in G\}
\label{eq:funalbe}\\
\mbox{Fun}_{\om} (G\times G_{AB},\Valbe) &:=& \{\phi:G\times G_{AB}\to
\Valbe \;|\; \phi(x\zinv, \xi) = \om(z) \phi(x,\xi), \nn\\
& &\mbox{for}\;z\in Z\}
\label{eq:funom}
\eea
We will also need the following sets:
\bea
G_o &:=& \{x\in G\,|\, \mbox{if}\; n_1\in N_A,\; n_2\in N_B\;\mbox{and}
\; n_1 x \ninv_2 = x \;\mbox{then}\; n_1 = n_2 \in Z\}
\label{eq:Gnot}\\
(G\times G)_o &:=& \{ (u,v)\in G\times G\,|\, u^{-1}v \in G_o\}\\
(G_{AB})_o &:=& \{ \xi\in G_{AB}\,|\, y(\xi)\in G_o\}
\eea
They have the following properties, which can be 
easily verified:\\
(a) If $x\in G_o, n_1\in N_A, n_2\in N_B$ then $n_1 x \ninv_2 \in G_o$. \\
(b) $x\in G_o \; \Leftrightarrow\; y(N_A x N_B)\in G_o$. \\  
(c) If $x\in G_o$, and $m_1, n_1 \in N_A, m_2, n_2\in N_B$ then
\be
m_1 x m_2^{-1} = n_1 x\ninv_2 \Rightarrow \exists z\in Z\;\mbox{such
that}\; m_1 = n_1 z, m_2 = n_2 z
\ee
(d) If $\xi\in(G_{AB})_o$ then $\exists z\in Z$ such that $n_1(y(\xi)] =
n_2(y(\xi)] = z$. \\
(e) If $(u,v)\in (G\times G)_o$, and $m_1\in N_A,\;m_2\in N_B$ then
$\exists z\in Z$ such that
\be
n_1(m_1\uv\minv_2) = z m_1 n_1(\uv),\;\; n_2(m_1\uv\minv_2) = z m_2 n_2(\uv)
\label{eq:nzm}
\ee

The next Proposition is the {\sl first step} in the tensor product 
decomposition.
Roughly speaking, we will consider the functions $\Phi$ in 
Eq.(\ref{eq:funalbe})
as elements of the tensor product representation space. After restriction to
$(G\times G)_o$ these functions $\Phi$ can be rewritten in a bijective
linear way as functions $\phi$ in Eq.(\ref{eq:funom}), restricted to
$G\times (G_{AB})_o$. The action of $\DG$ on $\Phi$ affects both arguments
of $\Phi$ (according to Eq.(\ref{eq:tensproduct})), but the corresponding
action on $\phi$ only affects its first argument, as we will see in 
Lemma \ref{lemtwo}. The second argument will in fact be directly related to 
the conjugacy class part of the label $(C(\xi),\gm)$ of a `new' irreducible 
representation of $\DG$, and thus we will prove that the tensor product
representation space is isomorphic to a direct integral of 
representation spaces of $\Pi^{C(\xi)}_{\om}$, where
$\Pi^{C(\xi)}_{\om}$ is not yet irreducible.

\begin{prop}\label{propmaps}
There is a linear map
\be
{\cal F}_1\,:\,\Phi \mapsto \phi\,:\, \mbox{Fun}_{\al,\bt}(G\times G,
\Valbe) \to \mbox{Fun}_{\om}(G\times G_{AB},\Valbe)
\label{eq:offmapone}
\ee
given by
\be
\phi(x,\xi) := \Phi(x w(\xi)^{-1}, xw(\xi)^{-1} y(\xi)),\quad x\in G,\;
 \xi \in G_{AB}.
\label{eq:mapone}
\ee
This map, when considered as a map 
\be
\Phi\mapsto\phi\,:\, 
\mbox{Fun}_{\al,\bt}((G\times G)_o,\Valbe) \to \mbox{Fun}_{\om}(G\times 
(G_{AB})_o,\Valbe),
\ee 
is a linear bijection with inversion formula ${\cal F}_2\,:\,\phi\mapsto
\Phi$ given by
\be
\Phi(u,v) = \al(n_1(\uv))\otimes \bt(n_2(\uv))\; \phi(un_1(\uv) w(N_A\uv N_B),
N_A \uv N_B).
\label{eq:maptwo}
\ee
\end{prop}
\Proof 
(i) Let $\phi$ be defined in terms of $\Phi\in \mbox{Fun}_{\al\bt}(G\times G,
\Valbe)$ by Eq.(\ref{eq:mapone}). The covariance condition of $\phi$ w.r.t.\
$Z$ follows because, for $z\in Z$,
\bea
\phi(x\zinv,\xi) &=& \Phi(x\zinv \wxinv, x\zinv \wxinv y(\xi)) = 
\Phi(x\wxinv \zinv, x \wxinv y(\xi) \zinv) = \nn\\
&=& \al(z)\otimes \bt(z)\, \Phi(x\wxinv, x\wxinv y(\xi)) = \om(z) \phi(x,\xi).
\nn
\eea
Moreover, $\phi$ restricted to $G\times (G_{AB})_o$ only involves $\Phi$
restricted to $(G\times G)_o$, since for $\xi\in(G_{AB})_o$ we have that 
$(x\wxinv, x\wxinv y(\xi))\in (G\times G)_o$.\\
(ii) ${\cal F}_1$ is injective because 
\bea
\left(({\cal F}_2\circ {\cal F}_1)\Phi\right)(u,v) &=&
\al(n_1(\uv)) \otimes \bt(n_2(\uv))\, \Phi(un_1(\uv), un_1(\uv) 
 y(N_A\uv N_B))=\nn\\
 &=&\Phi(un_1(\uv) n_1(\uv)^{-1}, un_1(\uv) y(N_A\uv N_B) n_2(\uv)^{-1})=\nn\\
 &=& \Phi(u, u \uv) = \Phi(u,v)\nn
\eea
and thus ${\cal F}_2\circ {\cal F}_1=$id. (Here it is not yet
necessary to restrict $(u,v)$ to $(G\times G)_o.)$ \\
(iii) Let $\Phi$ be defined in terms of $\phi \in\mbox{Fun}_{\om}(G\times 
(G_{AB})_o;\Valbe)$ by Eq.(\ref{eq:maptwo}). 
The covariance condition of $\Phi$ w.r.t.\ $N_A\times N_B$
follows because, for $m_1\in N_A, m_2\in N_B$ and $(u,v) \in (G\times
G)_o$ 
\bea
\Phi(um_1^{-1}, vm_2^{-1}) &=& \al(n_1(m_1\uv m_2^{-1}))\otimes \bt(n_2(m_1
\uv m_2^{-1})) \qquad\qquad\qquad\nn\\
& &\qquad\qquad\phi(u\minv_1 n_1(m_1 \uv \minv_2) w(N_A\uv N_B), 
N_A\uv N_B)\nn\\
&=& (\al(z)\otimes \bt(z))(\al(m_1)\otimes \bt(m_2)) (\al(n_1(\uv))\otimes
\bt(n_2(\uv))) \nn\\
& &\qquad \qquad\phi(uzn_1(\uv) w(N_A\uv N_B), N_A\uv N_B) =\nn\\
&=& (\al(m_1)\otimes \bt(m_2)) (\al(n_1(\uv))\otimes\bt(n_2(\uv)))\nn\\
& &\quad\phi(un_1(\uv) w(N_A\uv N_B), N_A\uv N_B) = 
\al(m_1)\otimes \bt(m_2)\, \Phi(u,v)\nn
\eea
for some $z\in Z$, where we have used property (e) from above. \\ 
(iv) ${\cal F}_1$ is surjective (or:  ${\cal F}_2$ is injective)
because for $(x,\xi)\in G\times (G_{AB})_o$ 
\bea
\left(({\cal F}_1\circ {\cal F}_2)\phi\right)(x,\xi) &=&
\al(n_1(y(\xi)))\otimes\bt(n_2(y(\xi)))\, \phi(x\wxinv n_1(y(\xi)) w(\xi), 
\xi) = \nn\\
&=& \al(z)\otimes \bt(z)\, \phi(xz,\xi) = \phi(x,\xi)\nn
\eea
for some $z\in Z$, where we have used property (d) from above. This 
concludes the proof.~\eof
\vspace*{.2cm}

Define a Borel measure $\mu$  such that 
\be
\int_G\,f(N_A y N_B)\,dy = \int_{G_{AB}}\, f(\xi)\,d\mu(\xi)
\label{eq:measure}
\ee
for all $f\in C(G_{AB})$. The measure $\mu$ has support $G_{AB}$.\\
We will now specialise the map ${\cal F}_1$ from Eq.(\ref{eq:offmapone}) to
the $L^2$-case.
\be
{\cal F}_1 \,:\,\Phi \mapsto \phi\,:\, L^2_{\al,\bt}(G\times G: \Valbe)
\to L^2_{\om}(G\times G_{AB},\Valbe).
\ee
Here the first $L^2$-space is defined as the representation space of a
tensor product representation, see Eqs.(\ref{eq:tenscov}) and
(\ref{eq:tensnorm}), and
the second $L^2$-space is defined as the set of all measurable
$\phi\,:\, G\times G_{AB} \to \Valbe$ satisfying, for all
$z\in Z$ that $\phi(x\zinv, \xi) = \om(z) \phi(x,\xi)$ almost
everywhere, and such that
\be
\|\phi\|^2 := \int_{\xi \in G_{AB}} \int_{x\in G} \|\phi(x,\xi)\|^2\,
dx\,d\mu(\xi)\,< \infty,
\ee
with almost equal $\phi$'s being identified. 
\begin{lem}\label{lemone}
Let $\Phi \in \mbox{Fun}_{\al,\bt}(G\times G,\Valbe)$
and let $\phi$ be given by Eq.(\ref{eq:mapone}). If $\Phi:G\times G\to 
\Valbe$ is moreover Borel measurable then $\phi:G\times G_{AB}\to \Valbe$ is
Borel measurable, and
\be
\int_{\xi\in G_{AB}}\int_{x\in G} \|\phi(x,\xi)\|^2\,dx\,d\mu(\xi) = 
 \int_G \int_G \|\Phi(u,v)\|^2\, du\,dv
\label{eq:phinorms}
\ee
In particular, the map ${\cal F}_1 :\Phi \mapsto \phi$ is an isometry of
the Hilbert space $L^2_{\al,\bt}(G\times G,\Valbe)$ into (not necessarily
onto!) the Hilbert space $L^2_{\om}(G\times G_{AB},\Valbe)$.
\end{lem}
\Proof
It follows from Eq.(\ref{eq:mapone}) and Proposition \ref{prop:maps}(a)
that $\phi$ is Borel measurable if $\Phi$ is Borel measurable. 
The left hand side of Eq.(\ref{eq:phinorms}) equals
\bea
\int_{G_{AB}}\left(\int_G \|\Phi(x\wxinv,x\wxinv y(\xi))\|^2\,
dx\right)\,d\mu (\xi) =
\int_G \int_{G_{AB}} \|\Phi(u, u y(\xi))\|^2\,du\,d\mu (\xi)=\nn\\
=\int_G \int_G \|\Phi(u,uy(N_A v N_B)\|^2 du\,dv = \int_G\left(\int_G
\|\Phi(u,un_1(v)^{-1}v n_2(v))\|^2 du\right)dv=\nn\\
=\int_G \int_G \|\Phi(u,uv)\|^2 \,du\,dv = \int_{G\times G} 
\|\Phi(u,v)\|^2\,du\,dv \qquad\nn ~\eof
\eea

Subsequently we can show how the map ${\cal F}_1$ transfers the action of
$\DG$ on $\Phi$ to an action of $\DG$ on $\phi$:
\begin{lem}\label{lemtwo}
Let $\Phi \in L^2_{\al,\bt}(G\times G,\Valbe), F\in\DG$ and 
\be
\Psi := (\PAal\otimes \PBbe)(F) \Phi. 
\ee
Let $\phi$ be defined in terms of $\Phi$ and
$\psi$ in terms of $\Psi$ via Eq.(\ref{eq:mapone}). Then
\be
\psi(x,\xi) = \int_G F(x g_{C(\xi)}\xinv, w)\,\phi(w^{-1}x,\xi)\,dw
\label{eq:psiphi}
\ee
\end{lem}
\Proof
\bea
\psi(x,\xi) &=& \left((\PAB)(F)\Phi\right)(x\wxinv, x\wxinv y(\xi))=\nn\\
&=&\int_G F(x\wxinv g_A y(\xi) g_B \yxinv w(\xi)\xinv, w)\,\Phi(w^{-1} x
\wxinv, w^{-1} x \wxinv y(\xi))\,dw = \nn\\
&=& \int_G F(xg_{C(\xi)}\xinv, w)\,\phi(w^{-1}x,\xi)\,dw\qquad\nn ~\eof
\eea

For $C\in\Conj(G)$ define a $*$-representation $\Pi^C_{\om}$ of $\DG$
on $L^2_{\om}(G,\Valbe)$ as follows:
\bea
\left(\Pi^C_{\om} (F)\phi\right)(x) &:=& \int_G F(x g_C\xinv,w)\,\phi(w^{-1}x)
\,dw,
\label{eq:piomega}\\
& &\qquad F\in\DG, \phi\in L^2_{\om}(G,\Valbe)\nn
\eea
This has the same structure as the defining formula for the representation
$\PAal$ as given in
Eq.(\ref{eq:representation}), but the covariance condition on the functions
$\phi$ in Eq.(\ref{eq:piomega}) is weaker, because it only involves right
multiplication of the argument with respect to $z\in Z$. 
Eq.(\ref{eq:psiphi}) can also be formulated as:
\be
\psi(x,\xi) = \left(\Pi^{C(\xi)}_{\om}(F) \phi(.,\xi)\right)(x)
\label{eq:psiphiprime}
\ee
which clearly shows that Lemmas \ref{lemone} and \ref{lemtwo} form the
first step in a direct integral decomposition of the representation
$\PAB$ into irreducible representations.
We will also need the following
\begin{ass} \label{ass:Go}
The complement of $G_o$ has measure zero in $G$.
\end{ass} 
This implies that
the complement of $(G\times G)_o$ has measure zero in $G\times G$, and the
complement of $(G_{AB})_o$ has measure zero in $G_{AB}$. 
For $G=SU(n)$ or $U(n)$ this Assumption will be satisfied if $A$ and
$B$ are conjugacy classes for which $g_A$ and $g_B$ are diagonal
matrices with all diagonal elements distinct (so they are regular
elements of the maximal torus $T$ consisting of diagonal
matrices). Then $N_A=N_B=T$, and $G_o$ certainly contains all
$g=(g_{ij}) \in G$ which have only nonzero off-diagonal elements, so
for which $g_{ij}\neq 0$ if $i\neq j$. Clearly, Assumption
\ref{ass:Go} is then satisfied. 
\begin{cor}\label{cor:onto}
The `isometry into' of Lemma \ref{lemone} can be narrowed down to an
`isometry onto', namely; \\
The map ${\cal F}_1\,:\,\Phi\mapsto \phi$ is an isometry of
the Hilbert space $L^2_{\al,\bt}((G\times G)_o;\Valbe)$ onto the
Hilbert space $L^2_{\om}(G\times (G_{AB})_o;\Valbe)$.
\end{cor}
\bigskip

The {\sl second step} 
in the decomposition of the tensor product representation $\PAB$  is the
decomposition of the representation $\Pi^{C(\xi)}_{\om}$ into irreducible 
components $\Pi^C_{\gm}$. In other words, to decompose the action of
$\DG$ on $L^2_{\om}(G\times G_{AB},\Valbe)$  as given by Eq.(\ref{eq:psiphi})
or Eq.(\ref{eq:psiphiprime}).
For the moment  suppose that $\xi$ can be fixed in 
Eq.(\ref{eq:psiphiprime}). Comparison of Eq.(\ref{eq:psiphiprime}) and
Eq.(\ref{eq:piomega}) with Eq.(\ref{eq:representation}) then shows
that essentially we have to decompose 
$L^2_{\om}(G)$ \footnote{The fact that the elements of $L^2_{\om}(G)$
should map to $\Valbe$ is not important for this argument}
as a direct sum of Hilbert spaces $L^2_{\gm}(G,V_{\gm})$  (possibly with
multiplicity) on which $\DG$ acts by the irreducible representation 
$\Pi^{C(\xi)}_{\gm}$, with $\gm \in \hat{N}_{C(\xi)}$.  

For $\phi\in L^2_{\om}(G)$, $C\in\Conj(G)$, $\gm\in\hat{N}_C$, $d_\gm:=
\dim V_\gm$, and $i,j=1,\ldots,d_\gm$ put
\be
\phi_{ij}^{C,\gm}(x):= \int_{N_C}\,\gm_{ij}(n)\,\phi(xn)\,dn, \quad x\in G,
\label{eq:psiij}
\ee
where we have chosen an orthonormal basis of $V_{\gm}$. 
By construction, for $n\in N_C$  we have that
\be
\phi^{C,\gm}_{ij}(xn) = \sum_{k=1}^{d_{\gm}} \gm_{ik}(\ninv) \phi^{C,\gm}_{kj}
(x).
\label{eq:phicovar}
\ee
For each $j=1,\ldots, d_{\gm}$ the vector $\phi^{C,\gm}_{ij}(x)$ 
takes values in $V_{\gm}$, the label $i$ denoting the component. Thus  
\be
\left(\phi^{C,\gm}_{ij}\right)_{i=1,\ldots,d_{\gm}} \in L^2_{\gm}(G,V_{\gm}).
\ee
Also, for $F\in \DG$
\be
\left(\left(\Pi^C_{\om}(F) \phi\right)^{C,\gm}_{ij}\right)_{i=1,\ldots,d_{\gm}}
= \Pi^C_{\gm}(F)\left(\phi^{C,\gm}_{ij}\right)_{i=1,\ldots,d_{\gm}},
\label{eq:actionphiij}
\ee
which follows from combining Eqs.(\ref{eq:representation}), (\ref{eq:piomega})
and (\ref{eq:psiij}). 

However, not all $\gm\in \hat{N}_C$ will occur, because $\phi^{C,\gm}_{ij}=0$
if $\gm|_Z \neq \om\,\mbox{id}$. This follows from the observation that
\bea
\phi^{C,\gm}_{ij}(x) &=& \int_{N_C} \gm_{ij}(n) \phi(xn)\,dn = \int_Z
\int_{N_C} \gm_{ij}(nz) \phi(xnz)\,dn\,dz=\nn\\
&=&\int_{N_C}\left(\int_Z \gm_{ij}(nz) \om(\zinv)\,dz\right)\,\phi(xn)\,dn = 
\left(\int_Z \gm(z)\om(\zinv)\,dz\right)\, \phi^{C,\gm}_{ij}(x).
\eea
Thus we must take $\gm$ to be an element of 
\be
\NCom =\{ \gm\in\hat{N}_C\;|\;\gm|_Z = \om\,\mbox{id}\}.
\label{eq:setNCom}
\ee

{}From the Peter-Weyl theorem applied to the function
$n\mapsto\phi(xn)$ with $x\in G$ we can derive that
\be
\int_G\, \|\phi (x)\|^2\,dx = \sum_{\gm\in\NCom} d_{\gm}
\sum_{i,j=1}^{d_{\gm}}\int_G\,\|\phi_{ij}^{C,\gm}(x)\|^2\,dx.
\label{eq:PW}
\ee
Thus as a continuation of the maps in Proposition \ref{propmaps}  we have an 
isometry 
\be
{\cal G}_1\,:\,\phi \mapsto \left(\sqrt{d_{\gm}} \left(\phi^{C,\gm}_{ij}
\right)_{i=1,\ldots,d_{\gm}} \right)_{\gm\in\NCom, j=1,\ldots,d_{\gm}}
\label{eq:G1}
\ee
of the Hilbert space $L^2_{\om}(G)$ into the direct sum of (degenerate)
Hilbert spaces
\be
\bigoplus_{\gm\in\NCom} \left(L^2_{\gm}(G,V_{\gm})\right)^{d_{\gm}}
\ee
which is intertwining between the representations $\Pi^C_{\om}$ and
$\oplus_{\gm\in\NCom} d_{\gm} \PCga$ of $\DG$. 

{}From the existence of an inversion formula we can see that the map 
${\cal G}_1$ is even an isometry {\it onto}. To that aim, fix $\gm \in \NCom$ 
and take $\left(\psi_{ij}\right)_{i,j=1,\ldots,d_{\gm}} \in\left(L^2_{\gm}
(G,V_{\gm})\right)^{d_{\gm}}$, i.e.\ $\psi_{ij} \in L^2(G)$ for $i,j=1,\ldots,
d_{\gm}$ and
\be
\psi_{ij}(xn) = \sum_{k=1}^{d_{\gm}} \gm_{ik}(\ninv) \psi_{kj}(x),\quad 
n\in N_C.
\ee
The map
\be
{\cal G}^{\gm}_2 \;:\;\psi\mapsto \phi\,:\,\left(L^2_{\gm}(G,V_{\gm})
\right)^{d_{\gm}}\,\to\,L^2_{\om}(G)
\ee
is defined by
\be
\phi(x) := d_{\gm} \sum_{k=1}^{d_{\gm}} \psi_{kk}(x).
\label{eq:G2}
\ee
Then indeed $\phi(x\zinv) = \om(z) \phi(x)$ with $z\in Z$.
Furthermore we have that ${\cal G}_1 \circ {\cal G}_2 =\,\mbox{id}$, since
for $\dl \in\left(\hat{N}_C\right)_{\om}$ and $\phi$ given by 
Eq.(\ref{eq:G2}) we have
\bea
\phi^{C,\dl}_{ij}(x) &=& d_{\gm}\sum_{k=1}^{d_{\gm}}\int_{N_C} \dl_{ij}(n)
\psi_{kk}(xn)\,dn \nn\\
&=& d_{\gm} \sum_{k,l=1}^{d_{\gm}} \left(\int_{N_C} \dl_{ij}(n) 
\overline{\gm_{lk}(n)}\,dn\right)\,\psi_{lk}(x) =
\left\{\begin{array}{cc}\psi_{ij}(x), & \dl = \gm \\ 0, & \dl\neq\gm.
\end{array}\right.
\eea

We want to apply the above decomposition of $L^2_{\om}(G)$ to our case
of $L^2_{\om}(G\times G_{AB},\Valbe)$. A slight problem occurs since 
in Eq.(\ref{eq:psiphiprime}) we had fixed $\xi$, which is not allowed
in an $L^2$-space. For varying $\xi$ we will have varying $C(\xi)$ and
hence varying $N_{C(\xi)}$ and $\left(\hat{N}_{C(\xi)}\right)_{\om}$. In
order to keep this under control we make the following 
\begin{ass}\label{asssplit}
$\Conj(G)$ splits as a disjoint union of finitely many Borel sets
$\Conj_p(G)$, on each of which $N_C$ does not vary with $C$. 
\end{ass}
For $G=SU(n)$ or $U(n)$ this assumption certainly holds, because we
can take the representatives $g_C =
\mbox{diag}(e^{i\th_1},...,e^{i\th_n})$ with
$\th_1\leq\th_2\leq...\leq \th_n <\th_1+2\pi$. Then $N_C$ only depends
on the partition of the set $\{1,...,n\}$ induced by the equalities
or inequalities between the $\th_j$'s. 

\input amssym.def
\def\al{\alpha}
\def\gog{{\goth g}}
\def\Cc{{\Bbb C}}

We would like to know whether the assumption holds for general
compact connected Lie groups $G$. Let $T$ be a maximal torus in $G$.
For any conjugacy class $A$ in $G$
take the representative $g_A$ uniquely as an element $t\in \overline{K}
\subset T$ (see after Assumption \ref{ass:conjcont}).
Van den Ban \cite{vdBan} has described the centraliser of $t$
in  $G$. From \cite{vdBan} we
conclude that the possible centraliser subgroups form a finite collection.
This can be seen as follows. Let $\gog_{\Cc}$ be the complexified
Lie algebra of $G$, let $\Sg$ be the root system of $T$ in $\gog_{\Cc}$,
and let $\gog_{\al}$ be the root space for $\al\in\Sg$.
Let $W$ be the Weyl group of the root system $\Sg$, which can also be
realized as the quotient group $W={\rm normaliser}_{G(T)}/T$. 
Let $t\in T$. Then the centraliser of $t$ in $G$ is completely
determined by the two sets (each a finite subset of a given finite set):
\be
\Sg(t):=\{\al\in\Sg \mid {\rm Ad}(t)\,X=X\;{\rm for}\,X\in\gog_{\al}\},
\quad
W(t):=\{w\in W\mid wtw^{-1}=t\}.
\ee
This also shows that, for $t_0\in T$, the set
$\{t\in T\mid \Sg(t)=\Sg(t_0),\;W(t)=W(t_0)\}$ is Borel.
Thus Assumption \ref{asssplit} is satisfied if $G$ is a compact connected  
Lie group. 
Note that the Lie algebra of the centralizer of $t$ in $G$ is 
determined by $\Sg(t)$ (see for instance Ch.~V, Proposition (2.3) in 
\cite{BrockDieck}).
For determining the centralizer itself, we need also $W(t)$.
This can be seen (cf.\ \cite{vdBan}) by using the so-called Bruhat
decomposition for a suitable complexification $G_{\Cc}$ of $G$.
\vspace*{.2cm}

Put $N_C = N_p$ if $C\in \Conj_p(G)$ and $G_{AB,p} :=\{\xi \in G_{AB}\,|\,
C(\xi)\in \Conj_p(G)\}$. Similarly to Eq.(\ref{eq:psiij}) for any
$\phi\in L^2_{\om}(G\times G_{AB},\Valbe)$ we define
\be
\phi^{p,\gm}_{ij}(x,\xi) := \int_{N_p}\gm_{ij}(n)\,\phi(xn,\xi)\,dn,\quad
x\in G, \xi\in G_{AB, p}, \gm\in\left(\hat{N}_p\right)_{\om},\, i,j=1,\ldots,
d_{\gm}.
\label{eq:phip}
\ee
with of course the same right covariance as Eq.(\ref{eq:phicovar}). 
Because $\phi$ now maps to $\Valbe$ we can say that
\be
\left(\phi^{p,\gm}_{ij}\right)_{i=1,\ldots,d_{\gm}} \in L^2_{\gm}(G\times
G_{AB, p},\Valbe\otimes V_{\gm})
\ee
where again $i$ denotes the component in $V_{\gm}$. Eq.(\ref{eq:actionphiij})
can now be generalised to
\be
\left(\left(\Pi^{C(\xi)}_{\gm}(F)\,\phi(.,\xi)\right)^{p,\gm}_{ij}\right)_{i=
1,\ldots,d_{\gm}} = \Pi^{C(\xi)}_{\gm} (F)\left(\phi^{p,\gm}_{ij}(.,\xi)
\right)_{i=1,\ldots,d_{\gm}},\quad \xi\in G_{AB, p}.
\ee

Corresponding to Eq.(\ref{eq:PW}) we now have the isometry property
\be
\int_G \int_{G_{AB}}\|\phi(x,\xi)\|^2\,dx\,d\mu(\xi) = 
\sum_p \sum_{\gm\in\left(\hat{N}_p\right)_{\om}} d_{\gm} \sum_{i,j=1}^{d_{\gm}}
\int_G \int_{G_{AB,p}}\|\phi^{p,\gm}_{ij}(x,\xi)\|^2\,dx\,d\mu(\xi)
\label{eq:phipnorm}
\ee
and the isometry from Eq.(\ref{eq:G1}) now becomes the isometry
\be
{\cal G}_1\,:\,\phi \mapsto \left(\sqrt{d_{\gm}} \left(\phi^{p,\gm}_{ij}
\right)_{i=1,\ldots,d_{\gm}} \right)_{p;\gm\in (\hat{N}_p)_{\om}, j=1,\ldots,
d_{\gm}}
\label{eq:G1p}
\ee
of the Hilbert space $L^2_{\om}(G\times G_{AB},\Valbe)$ into the
direct sum of Hilbert spaces 
\be
\bigoplus_p \bigoplus_{\gm\in (\hat{N}_p)_{\om}}
\left(L^2_{\gm}(G\times  G_{AB,p},\Valbe\otimes V_{\gm})\right)^{d_{\gm}}.
\label{eq:sumhilbsp}
\ee
This isometry is intertwining between the direct integral of representations 
\be
\int_{G_{AB}}^{\oplus}\Pi^{C(\xi)}_{\om} d\mu(\xi)\qquad
\mbox{and}\qquad
\bigoplus_p \bigoplus_{\gm\in (\hat{N}_p)_{\om}}
\int_{G_{AB,p}}^{\oplus}d_{\al} d_{\bt} d_{\gm}\,\Pi^{C(\xi)}_{\gm}\,d\mu(\xi)
\ee
of $\DG$. Keep in mind that only the equivalence class of the measure 
$\mu$ matters in a direct integral of representations, as above. 

Again, to show that ${\cal G}_1$ is indeed an isometry into, we
construct the inverse: \\
for $(\psi_{ij})_{i,j=1,\ldots,d_{\gm}}\in L^2_{\gm}(G\times G_{AB, p}:
\Valbe\otimes V_{\gm})^{d_{\gm}}$ define the map
\be
{\cal G}^{p,\gm}_2\;:\;\psi \mapsto \phi\;:\;L^2_{\gm}(G\times G_{AB, p}:
\Valbe\otimes V_{\gm})^{d_{\gm}} \to L^2_{\om}(G\times G_{AB},\Valbe)
\ee
by
\be
\phi(x,\xi) := d_{\gm}\sum_{k=1}^{d_{\gm}} \psi_{kk}(x,\xi).
\ee
Then ${\cal G}_1 \circ {\cal G}^{p,\gm}_2 =\,\mbox{id}$, which can be 
shown in the same way as under Eq.(\ref{eq:G2}).
\vspace*{.2cm}

We now combine {\sl step one} and {\sl step two} in the procedure 
described above. The decomposition of the tensor product representation is
then given by the intertwining isometry $ \rho := {\cal G}_1 \circ
{\cal F}_1$, and its  
inverse is given by ${\cal F}_2 \circ {\cal G}^{p,\gm}_2$. (The latter acting 
on $L^2_{\gm} (G\times G_{AB, p},\Valbe\otimes V_{\gm})^{d_{\gm}}$.)
 
Thus we have determinded the Clebsch--Gordan series from
Eq.(\ref{eq:fusionintegral})  
\be
\PAB \simeq  \int_{G_{AB}}^{\oplus}
 \,\bigoplus_{\gm\in\hat{N}_C}\, d_{\al} d_{\bt} d_{\gm}
 \Pi^{C(\xi)}_{\gm}\,d\mu(\xi),
\label{eq:seriesgeneral} 
\ee
with $\mu$ an equivalence class of measures. More precisely, we have
to take the variation of $N_{C(\xi)}$ with $\xi$ into account, which
splits the direct integral over $\xi$:
\be
\PAB \simeq  \bigoplus_p \bigoplus_{\gm\in (\hat{N}_p)_{\om}}
\int_{G_{AB,p}}^{\oplus}d_{\al} d_{\bt} d_{\gm}\,\Pi^{C(\xi)}_{\gm}\,d\mu(\xi)
\ee 
 
Combining ${\cal F}_1$ from Eq.(\ref{eq:mapone}) and ${\cal G}_1$ from 
Eq.(\ref{eq:G1p}) we see that a  $\Phi \in L^2_{\al\bt}(G\times G:
\Valbe)$ is taken to an `object' in the direct sum/integral of Hilbert spaces 
\be
\bigoplus_p \bigoplus_{\gm\in(\hat{N}_p)_{\om}}\int^{\oplus}_{G_{AB,p}} 
L^2_{\gm}(G\times G_{AB,p},\Valbe\otimes V_{\gm})^{d_{\gm}}\,d\mu(\xi)
\ee
This object  depends on 
$\xi \in G_{AB}$, which determines the class label $C$ of the (irreducible)
representation $\Pi^{C(\xi)}_{\gm}$ which occurs in the decomposition. It 
has an index $i$ denoting the component of the vector (with tensor products
of vectors in $\Valbe$ as its entries) in $V_{\gm}$ to 
which a group element $x$ is mapped, an index $p$ which
denotes the Borel set in $\Conj(G)$, which in turn determines
the set $(\hat{N}_p)_{\om}$ to which the label $\gm$ of the
$\DG$-representation 
must belong. Finally, the object has an index $j$ indicating the degeneracy 
of the irreducible representation $\Pi^{C(\xi)}_{\gm}$. The `vector of tensor
products of vectors' means that each component in $V_{\gm}$ of the object 
in fact depends on the full vector in $\Valbe$ to which $\Phi$ maps a pair
$(x_1, x_2)\in G\times G$. We can `dissect' the isometry $\rho$ according to
the way it maps the components of $\Phi$ to components of the object described
above, this results in the following
\begin{thm}\label{thm:rho}
Let $\PAal, \PBbe$  be irreducible $*$-representations of $\DG$, and let $p$
label the finitely many Borel sets in $\Conj(G)$, on each of which $N_C$ does
not vary with $C$. Take $\xi\in G_{AB, p}$ and $\gm \in
(\hat{N}_p)_{\om}$. Then, for each  $k=1,..., d_{\al}$ and $l=1,...,d_{\bt}$ 
and $i,j=1,...,d_{\gm}$ a mapping
\be
\rho^{\xi}_{\gm, k,l,j}\,:\,L^2_{\al\bt}(G\times G,V_{\al}\otimes V_{\bt})
\to L^2_{\gm}(G,V_{\gm}) 
\ee
intertwining the representations $\PAB$ and $\Pi^{C(\xi)}_{\gm}$ is given by
\bea
\left(\rho^{\xi}_{\gm,k,l,j}\Phi\right)_i (x) &:=& \left(
\phi^{p,\gm}_{ij}(x,\xi)\right)_{k,l} \nn\\
&=& \int_{N_{C(\xi)}} \gm_{ij}(n)\;\Phi_{kl}(xn\wxinv,xn\wxinv y(\xi))\,dn.
\label{eq:rhoxi}
\eea
\end{thm}
An implicit expression for the fusion rules (multiplicities) can now also 
be obtained by comparing the squared
norms before and after the action of $\rho$ on $\Phi$. We then would like to 
rewrite a direct integral over $G_{AB, p}$ of representations 
$\Pi^{C(\xi)}_{\gm}$  as a direct integral over $\Conj_p(G)$ of representations
$\PCga$. However, if the map $\xi\mapsto C(\xi)\,:\,G_{AB} \to \Conj(G)$ is
non-injective, which might be the case as we have mentioned before, this 
rewriting can be difficult. To solve this, 
we also define a Borel measure $\nu$ on $\Conj(G)$ such that
\be
\int_{G_{AB}} F(C(\xi))\,d\mu(\xi)=\int_{\Conj(G)}F(C)\,d\nu(C)
\label{eq:borelmeasure}
\ee
for all $F\in C(\Conj(G))$. The measure $\nu$ has support $\lmAB(G_{AB})$.
By Theorem \ref{113} there exists for almost each $C\in\Conj(G)$
a Borel measure $p^C$ on $G_{AB}$ such that
\be
\int_{G_{AB}} f(\xi)\,d\mu(\xi)=
\int_{C\in\Conj(G)}\left(\int_{\xi\in G_{AB}} f(\xi)\,
dp^C(\xi)\right)\,d\nu(C)
\label{118}
\ee
for each $f\in C(G_{AB})$. If the mapping $\lm_{AB}$ is 
injective (like in the case of $G=SU(2)$, as we will discuss in the
next section) then the above simplifies to
\be
\int_{G_{AB}} f(\xi)\,d\mu(\xi) = \int_{I_{AB}} f(\lm^{-1}_{AB} (C))\,d\nu(C),
\label{eq:measurenu}
\ee
where $I_{AB}$ is the image of $G_{AB}$ under $\lm_{AB}$.

Combining Eqs.(\ref{eq:phinorms}) and (\ref{eq:phipnorm}) the isometry
property which contains the implicit expression for the multiplicities
now reads
\bea\label{117}
\lefteqn{\int_G \int_G\|\Phi(u,v)\|^2\,du dv =\qquad} \\
& &\!\!\!\!\!d_{\al} d_{\bt}\sum_p \sum_{k=1}^{d_{\al}}
\sum_{l=1}^{d_{\bt}}  \sum_{\gm\in (\hat{N}_p)_{\om}}\!\!d_{\gm}
\sum_{j=1}^{d_{\gm}}\int_{\Conj_p(G)} 
\!\left( \int_{G_{AB, p}}\!\left(\sum_{i=1}^{d_{\gm}}
\int_G\|\left(\rho^{\xi}_{\gm, 
k,l,j}\,\Phi\right)_i(y)\|^2\,dy\right)dp^C(\xi)\right)d\nu (C).\nn
\eea
Eq.(\ref{117}) can be written more compactly as:
\be
\|\Phi\|^2 = d_{\al} d_{\bt} \sum_p\sum_{k=1}^{d_{\al}} 
\sum_{l=1}^{d_{\bt}} \int_{\Conj_p(G)}
\left(\sum_{\gm\in (\hat N_p)_{\om}}d_\gm\sum_{j=1}^{d_\gm}
\int_{G_{AB,p}}\|\rho^{\xi}_{\gm,k,l,j}\,\Phi\|^2\,dp^C(\xi)\right)d\nu(C).
\label{116}
\ee
If $\lm_{AB}$ is injective then Eq.(\ref{116}) simplifies to
\be
\|\Phi\|^2 = d_{\al} d_{\bt}\sum_p \sum_{k=1}^{d_{\al}}\sum_{l=1}^{d_{\bt}} 
\int_{I_{AB,p}}\left(
\sum_{\gm\in (\hat{N}_p)_{\om}}\,d_\gm\sum_{j=1}^{d_{\gm}}
\|\rho^{\lm_{AB}^{-1} 
(C)}_{\gm,k,l,j}\Phi\|^2\right)\,d\nu(C)
\label{eq:rhoCsimple}
\ee
with $I_{AB,p}=\lm_{AB}(G_{AB,p})$. Note that the measures no longer stand
for equivalence classes of measures, but for specific measures, since we
are comparing (squared norms of) vectors in Hilbert spaces. The measure
$\nu$ may involve a nontrivial Jacobian from the mapping $\lm_{AB}$. 

The multiplicities $\Nabc$ can now more or less be extracted from
Eq.(\ref{117}) or Eq.(\ref{116}), that is, we can conclude the following:\\
(i) $\Nabc=0$ if $C\not\in\lm_{AB}(G_{AB})$ \\
(ii) $\Nabc=0$ if $\gm\not\in \hat{N}_{\om}$ \\
(iii) if $\Nabc\neq 0$ then $\Nabc = d_{\al} d_{\bt} d_{\gm}$\\
(iv) the inner product on $V^C_{\gm}$ will depend nontrivially on $A$ and $B$
according to the Jacobian of the mapping $\lm_{AB}$ and its non-injectivity,
which is reflected in the measure $p^C(\xi)$. 
\section{Explicit results for $G = SU(2)$\label{s7}}
\setcounter{equation}{0}
To illustrate the above aspects of tensor products of irreducible 
representations we will now consider the case of $G = SU(2)$. We will
only discuss the decomposition of a `generic' tensor product representation
and give explicit formulas for the Clebsch--Gordan coefficients in this
case. Some applications and the treatment of more special tensor
products will be discussed elsewhere \cite{BM}.

In \cite{KM}
we have given the classification of the irreducible unitary representations
of ${\cal D}(SU(2))$. For application of the main result of this paper
(the decomposition of the tensor product of such representations
into single representations) we first need to establish the notation and
parametrisation of elements of $SU(2)$. In this section we  
use the conventions of Vilenkin \cite{Vil}, because this book contains
a complete and explicit list of formulas which are needed in our
analysis. For the Wigner functions we use
the notation of Varshalovich {\em et al} \cite{VMK} (especially chapter 4). 
\subsection{Parametrisation and notation}
To specify an $SU(2)$-element we use both the Euler angles
$(\phi,\th,\psi)$, and the parametrisation by a single rotation
angle $r$ around a given axis $\hat{n}$. In the Euler--angle 
parametrisation each $g\in SU(2)$ can be written as
\be
g = g_{\phi} a_{\th} g_{\psi}
\label{eq:euler}
\ee
with
\be
g_{\phi} = \left(\begin{array}{cc} e^{\half i \phi} & 0 \\ 0 &
e^{-\half i \phi}\end{array}\right) , \qquad
a_{\th} = \left(\begin{array}{cc} \cos \half\theta & -\sin \half\theta \\ 
\sin \half\theta & \cos \half\theta \end{array}\right) \nn
\label{eq:atheta}
\ee
\be
0\leq\th\leq \pi,\quad 0\leq\phi<2\pi,\quad -2\pi\leq\psi\leq 2\pi.
\label{eq:eulermatrices}
\ee
The diagonal subgroup
consists of all elements $g_{\phi}$, and is isomorphic to $U(1)$.

The conjugacy classes of $SU(2)$ are denoted by $C_r$ with $0\leq r\leq
2\pi$. The representative of $C_r$ can be taken to be $g_r$, so in the
diagonal subgroup. Then Assumption \ref{ass:conjcont} which states that
the map of the set of conjugacy classes of $G$ to $G$ itself (i.e.\ the
map to representatives) can be chosen to be continuous
is satisfied. For $r=0$ and $2\pi$ the centralizer
$N_0=N_{2\pi}=SU(2)$, 
for the other conjugacy classes the centralizer $N_r = U(1)$.

Let $0<r<2\pi$. Then $C_r$ clearly consists of the elements
\be
g(r,\th,\phi) := g_{\phi} a_{\th} g_r a_{\th}^{-1} g_{\phi}^{-1}. 
\label{eq:gagag}
\ee
If we take the generators of $SU(2)$ in the fundamental representation
to be
\bea
\tau_1:= \left(\begin{array}{cc}1&0\\0&-1\end{array}\right),\quad
\tau_2:= \left(\begin{array}{cc}0&1\\1&0\end{array}\right),\quad
\tau_3:= \left(\begin{array}{cc}0&i\\-i&0\end{array}\right)
\label{eq:tau}
\eea
and define the unit vector
\be
\hat{n}(\th,\phi) := (\cos\th, \sin\th \cos\phi, \sin\th \sin\phi)
\label{eq:direction}
\ee
then we can also write the element $g(r,\th,\phi)$ as
\be
g(r,\th,\phi) = \exp(i\frac{r}{2} \hat{n}(\th,\phi)\cdot\vec{\tau})=
\unit \cos\frac{r}{2} + \hat{n}\cdot\vec{\tau}\, i \sin\frac{r}{2}
\label{eq:exponent}
\ee

This means that there is a 1--1 correspondence between $\hat{n}(\th,\phi)$
and the cosets $g_{\phi} a_{\th} N_r$. In other words, the mapping
$\hat{n}(\th,\phi) \mapsto g(r,\th,\phi)\,:\, S^2 \to C_r$ is bijective
from the unit sphere $S^2$ in $\dubbelR^3$ onto the conjugacy class $C_r$.
\subsection{Irreducible representations}
Next we consider the `generic' irreducible unitary representations of 
${\cal D}(SU(2))$, i.e.\ for the case $r\neq 0,2\pi$. The other cases
will be treated
elsewhere \cite{BM}. The centralizer representations will
be denoted by $n\in\half\dubbelZ$ (so not the elements themselves as we 
did in the sections before, when we discussed the general case). The
irreducible unitary representations of $N_r$ are the 1-dimensional 
representations
\be
n\,:\,g_{\zeta}\mapsto e^{in\zeta} ,\qquad -2\pi\leq\zeta\leq 2\pi,\quad
n\in\half\dubbelZ.
\ee
For the generic representations $\Pi^r_n$ of ${\cal D}(SU(2))$ the 
representation space is
\be
V^r_n = \{\phi \in L^2(SU(2),\dubbelR/2\pi) \;|\;  
\phi(g g_{\zeta}) = e^{-in\zeta}\phi(g), \quad -2\pi\leq \zeta\leq 2\pi\}.
\label{eq:Vrn}
\ee
An orthogonal basis for $V^r_n$ is given by the Wigner functions $\D{j}{mn}$,
{\sl where the label $n$ is fixed}. A thorough treatment of the Wigner
functions as a basis of functions on $SU(2)$ can be found in \cite{Vil}.
For $g\in SU(2)$ parametrised by the Euler angles as in Eq.(\ref{eq:euler})
the Wigner function $D^j_{mn}$ corresponding to the $m,n$-th matrix element
in the $j$-th irreducible representation takes the value
\be
D^j_{mn}(g) = e^{-im\phi} P^j_{mn}(\cos \th) e^{-in\psi},
\ee
where $P^j_{mn}$ can be expressed in terms of Jacobi polynomials. For
all $g_{\zeta}=e^{i\zeta} \in U(1)$ we have that
\be
D^j_{mn} (xg_{\zeta}) = e^{-in\zeta} D^j_{mn}(x).
\ee
This shows indeed that the set
$\{D^j_{mn}\,| n\,\,\mbox{fixed}, j\in \half \dubbelN, j\geq n, -j\leq
m\leq j\}$ has the right covariance property.
The Wigner functions form a complete set on $SU(2)$, so the aforementioned
set forms a basis for a Hilbert space corresponding to an irreducible unitary
representation of ${\cal D}(SU(2))$, with fixed centraliser representation $n$
and arbitrary conjugacy class $0<r<2\pi$. In other words, the Hilbert spaces
for irreducible unitary representations with the same $n$ and different $r$
are equivalent, and thus can be spanned by identical bases. 
Recall that the $r$-dependence of the representation functions $\phi \in 
V^r_n$ is only reflected in the action of ${\cal D}(SU(2))$ on $V^r_n$:
\be
\left(\Pi^r_n (F)\phi\right)(y) = \int_{SU(2)} F(yg_r\yinv,x) \phi(\xinv y)\,
dx,\quad \phi \in V^r_n,
\ee
Strictly speaking, we should label the (basis) vectors
of $V^r_n$ by $r$ as well, then an arbitrary state in a generic
representation is written as
\be
\rphi_n (x)=\sum_{j>n}\sum_{-j\leq m\leq j} c_{jm}\,\rD^j_{mn}(x),
\quad x\in G.
\ee
(Note that the sum over $j$ is infinite.) However, since we will always 
specify which representation $\Pi^r_n$ we are dealing with, we will omit 
the $r$-label on the functions. 

By Eq.(\ref{eq:char}) the character $\chi^r_n$ of a generic representation
$\Pi^r_n$ is given by
\be
\chi^r_n(F) = \int_{SU(2)} \int_{U(1)} F(zg_r\zinv, zg_{\zeta}\zinv)\,
e^{in\zeta} d\zeta\,dz,\qquad F\in C^{\infty}(SU(2)\times SU(2)). 
\ee
\subsection{Clebsch--Gordan series}{\label{ssseries}}
First we will determine the decomposition of the tensor product of two
generic representations $\Pi^{r_1}_{n_1}$ and $\Pi^{r_2}_{n_2}$ as in
Eq.(\ref{eq:seriesgeneral}).  It will turn out that $p$ takes only one
value, corresponding to generic $r_3$, and that the map
$\lm_{r_1,r_2}$ is injective. We have to determine the image
$I_{r_1,r_2}$ of $\lm_{r_1,r_2}$, the equivalence class of the measure
$\nu$, and the set $\left(\hat{N}_{r_3}\right)_{\ep}$. Since the
centraliser representations $n_1, n_2, n_3$ are one-dimensional we see
that the nonvanishing multiplicities $N^{r_1 r_2 n_3}_{n_1 n_2 r_3}=1$. 

We choose $y(\th):= a_{\th}$ as a 
representative for the double coset $N_{r_1} a_{\th} N_{r_2}$, which
is an element of $G_{r_1 r_2}= U(1)\bsl SU(2)/U(1)$. Then 
Assumption \ref{ass2}, stating that the representatives of the double cosets
can be chosen in a continuous way, is satisfied. Eq.(\ref{eq:xyz}),
which for this case determines $r_3(\th)$ and $w(\th)$, now reads
\be
g_{r_1} a_{\th} g_{r_2} a^{-1}_{\th} = w(\th) g_{r_3(\th)} \winv(\th).
\label{eq:ys}
\ee
By computing the trace of the left-hand side of Eq.(\ref{eq:ys}) we 
find for $r_3= r_3(\th)$ that
\be
\cos\frac{r_3}{2} = \cos \frac{r_1}{2}\cos \frac{r_2}{2} - \cos\th
\sin \frac{r_1}{2} \sin \frac{r_2}{2},
\label{eq:rs}
\ee
which gives us the mapping $\lm_{r_1,r_2}$ from Eq.(\ref{eq:cosettoclass}):
\be
\lm_{r_1,r_2} (U(1) a_{\th} U(1)) =
2 \arccos (\cos\half r_1 \cos\half r_2 - \cos\th \sin\half r_1 \sin\half r_2).
\label{eq:lmr1r2}
\ee
Thus the mapping $\lm_{r_1,r_2}\,:\,G_{r_1 r_2}\,\to\, [0,2\pi]$
is injective with image 
\be
I_{r_1 r_2} = [|r_1 - r_2|, \mbox{min}(r_1+r_2, 4\pi-(r_1+r_2))].
\label{eq:imager1r2}
\ee 

Now we compute the measures $\mu$ and $\nu$ from
Eqs.(\ref{eq:borelmeasure}) and (\ref{118}). 
The measure $\mu$ on $G_{r_1 r_2}$ follows from
\be
\int_{SU(2)} f(g)\,dg = \half\int_0^{\pi} f(a_{\th}) \sin\th\,d\th
\label{eq:dmu}
\ee
for a function $f\in C(G_{r_1 r_2})$, and thus
\be
d\mu(\th) = \half \sin\th \,d\th.
\ee
The Borel measure $\nu$ on the set of conjugacy classes can be 
derived via
\be
\int_0^{\pi} F(\lm_{r_1,r_2}(U(1)a_{\th}U(1)))\,d\mu(\th) = \int_{I_{r_1,
r_2}} F(r_3)\,d\nu(r_3)
\ee
for an $F\in C(\Conj(SU(2)))$. With formula (\ref{eq:lmr1r2}) it follows that
\bea
d\nu (r_3)=\left\{\begin{array}{cl} 
\frac{\sin\frac{r_3}{2}}{4\sin\frac{r_1}{2} \sin \frac{r_2}{2}} dr_3 , &
|r_1-r_2|\leq r_3\leq\mbox{min}(r_1+r_2,4\pi-(r_1+r_2))\\ 0 , &
\mbox{otherwise}.\end{array}\right.
\label{eq:nur3}
\eea
We conclude that the nongeneric conjugacy classes $r_3=0$ and $r_3=2\pi$
have $\nu$-measure zero in $I_{r_1,r_2}$. We also see that
the measure $d\nu(r_3)$ is equivalent with the measure $dr_3$ on
$I_{r_1,r_2}$.

To determine $\left(\hat{N}_{r_3}\right)_{\ep}$ we remark that
\be
n_1(z)\otimes n_2(z) = \ep(z)\,\mbox{id}_{V_{n_1}\otimes V_{n_2}}, \quad
z = \{e, -e\}\subset SU(2).
\ee
So $\left(\hat{N}_{r_3}\right)_{\ep} = (n_1+n_2)\,\mbox{mod}\,\dubbelZ$.
The Clebsch--Gordan series now reads
\be
\Pi^{r_1}_{n_1}\otimes\Pi^{r_2}_{n_2} \simeq \bigoplus_{n_3\in
(n_1+n_2)\mbox{mod} 
\dubbelZ} \int^{\oplus}_{I_{r_1,r_2}} \Pi^{r_3}_{n_3}\, dr_3.
\label{eq:suseries}
\ee
\subsection{Clebsch--Gordan coefficients}
We will now explicitly construct the mapping $\rho$ from Eq.(\ref{eq:rhoxi}),
successively applying the steps of section \ref{s6}. We can
compute $w(\th)=g_{\phi_w} a_{\th_w}$ by first rewriting
Eq.(\ref{eq:ys}) as
\be
(\unit\cos\frac{r_1}{2} + i \tau_1\,
\sin\frac{r_1}{2})(\unit\cos\frac{r_2}{2}  + i (\cos\th\,\tau_1 +
\sin\th\, \tau_2)\sin\frac{r_2}{2})= \unit\cos\frac{r_3}{2}+
i\,\hat{n}_w\cdot\vec{\tau} \sin\frac{r_3}{2}
\ee
(in view of Eqs.(\ref{eq:gagag}), (\ref{eq:direction}),
(\ref{eq:exponent})), and then comparing coefficients of $\tau_1,
\tau_2, \tau_3$ on both sides. This yields
\bea
\hat{n}_{w(\th)} = \left(\begin{array}{c} \cos \th_w \\ \sin \th_w \cos\phi_w\\
\sin \th_w \sin\phi_w\end{array}\right) = 
\frac{1}{\sin\frac{r_3}{2}}  \left(
\begin{array}{c} \sin\frac{r_1}{2}\cos\frac{r_2}{2} +\cos\th\cos\frac{r_1}{2}
\sin\frac{r_2}{2}\\ \sin\th\cos\frac{r_1}{2} \sin\frac{r_2}{2}\\ \sin\th
\sin\frac{r_1}{2}\sin\frac{r_2}{2} \end{array}\right).
\label{eq:nw}
\eea
It follows from Eqs.(\ref{eq:ys}) and (\ref{eq:nw}) that $\th_w$ and
$\phi_w$ depend continuously on $\th$, even for $r_1=r_2$, in which
case the right hand side of Eq.(\ref{eq:nw}) tends to
\bea
\left(\begin{array}{c} 0 \\ \cos\frac{r_1}{2}\\
\sin\frac{r_1}{2}\end{array}\right)
\eea
as $\th \uparrow \pi$, hence $\th_w \to \frac{\pi}{2}, \phi_w \to
\frac{r_1}{2}$. Thus the Borel map from Proposition \ref{prop:maps}
(a) can be chosen continuously. 

The first step in the tensor product decomposition is the construction
of the map ${\cal F}_1$ from Corollary \ref{cor:onto}. The isometry 
\be
{\cal F}_1\;:\; L^2_{n_1,n_2}(SU(2)\times SU(2))\to L^2_{\ep}(SU(2)\times 
[0,\pi]),\qquad \ep = (n_1+n_2)\,\mbox{mod}\,\dubbelZ
\ee
is given by
\be
\phi(x,\th) = \Phi(x w(\th)^{-1}, x w(\th)^{-1} a_{\th}).
\ee
For the inversion formula ${\cal F}_2$ we need a choice for the Borel map 
from Proposition \ref{prop:maps} (b). It follows
straightforwardly from the Euler angle parametrisation:
write $x\in SU(2)$ as $x=g_{\phi_x} a_{\th_x} g_{\psi_x}S$ with
$0\leq \th_x\leq \pi,\;0\leq \phi_x< 2\pi,\; -2\pi\leq\psi_x <2\pi$. 
Put $y(U(1) x U(1)) := a_{\th_x}$ and $n_1(x):=g_{\phi_x},
\;n_2(x):=g_{\psi_x}$. Then ${\cal F}_2\,:\,\phi \mapsto \Phi$ is given by
\be
\Phi(u,v) = e^{in_1\phi_{\uv}} e^{in_2\psi_{\uv}}\,\phi(u g_{\phi_{\uv}}
w(\th_{\uv}), \th_{\uv})
\ee
with $\uv \in SU(2)_o$, and
\be
SU(2)_o = \left\{\left(\begin{array}{cc}\al & -\overline{\bt}\\ 
\bt&\overline{\al}\end{array}\right) \in SU(2)\;|\;\al,\bt\neq 0\right\}.
\label{eq:Gosu2}
\ee
Assumption \ref{ass:Go}, stating that the complement of $G_o$ has
measure zero in $G$, is satisfied for this case. 
\vspace*{.2cm}\\
The second step in the tensor product decomposition is given by the
isometry ${\cal G}_1$ from Eq.(\ref{eq:G1})
\be
{\cal G}_1\;:\;L^2_{\ep}(SU(2)\times [0,\pi]) \to \bigoplus_{n_3\in (n_1+n_2)
\,\mbox{mod}\,\dubbelZ} L^2_{n_3}(SU(2)\times I_{r_1, r_2}).
\ee
Assumption \ref{asssplit} about $\Conj(SU(2))$ is satisfied, because there are
two sets in $\Conj(SU(2))$ with distinct centralisers: the set $p_0 =
\{r=0, r=2\pi\}=Z$ with centraliser $SU(2)$, 
and the set $p_1 = \{r \in (0,2\pi)\}$ with centraliser $U(1)$. 
{}From Eq.(\ref{eq:nur3}) we see that the set $p_0$ will give no contribution
in the decomposition of the squared norm of the tensor product state, 
because for $r_3=0,2\pi$ the measure $\nu(r_3)$ on the conjugacy classes
is zero. Therefore we only need to compute Eq.(\ref{eq:phip}) for $p=p_1$:
\be
\phi^{p_1, n_3}(x,\th) = \int_{U(1)} e^{in_3\zeta} \phi (x g_{\zeta},
\th)\,d\zeta,\qquad n_3 \in (n_1+n_2)\;\mbox{mod}\;\dubbelZ
\ee
with the $U(1)$ over which we integrate embedded in $SU(2)$, so 
$-2\pi \leq \zeta \leq 2\pi$, and the Haar measure $d\zeta$ appropriately
normalised. 
The isometry property of Eq.(\ref{eq:phipnorm}) now becomes
\be
\int_{SU(2)} \int_0^{\pi}|\phi(x,\th)|^2\,dx\,d\mu(\th) =
 \sum_{n_3 \in (n_1+n_2)\mbox{mod}\dubbelZ} \int_{I_{r_1, r_2}}
\int_{SU(2)} |\phi^{p_1,n_3}(x, r_3)|^2\,dx\,d\nu(r_3).
\ee
The inverse mapping ${\cal G}^{p_1}_2$ reads
\be
\phi(x,\th) = \sum_{n_3\in (n_1+n_2)\;\mbox{mod}\;\dubbelZ} \phi^{p_1, 
n_3}(x,\th).
\ee
This results in the mapping $\rho$ intertwining the representations 
\be
\Pi^{r_1}_{n_1}\otimes \Pi^{r_2}_{n_2} \quad \mbox{and}\quad
\bigoplus_{n_3 \in (n_1+n_2)\;\mbox{mod}\;\dubbelZ} \int^{\oplus}_{I_{r_1,r_2}}
\Pi^{r_3}_{n_3}\,d\nu(r_3)
\ee
We calculate the components of mapping $\rho$ as given in Eq.(\ref{eq:rhoxi}).
The labels $i,j,k,l$ can 
be ignored, because $V_{n_1}, V_{n_2}, V_{n_3}$ are one-dimensional.
\be
\left(\rho^{\th}_{n_3}\Phi\right)(x) = \int_{U(1)} e^{in_3\zeta} 
\Phi(x g_{\zeta} w(\th)^{-1}, x g_{\zeta} w(\th)^{-1} a_{\th})\,d\zeta.
\label{eq:rhon3}
\ee
The Clebsch--Gordan series from Eq.(\ref{eq:suseries}) is contained in
\be
\int_{SU(2)}\int_{SU(2)} \|\Phi(u,v)\|^2\,du \,dv = \sum_{n_3 \in (n_1+n_2)\;
\mbox{mod}\;\dubbelZ}\int_{I_{r_1,r_2}}\left(\int_{SU(2)} |(\rho^{r_3}_{n_3}
\Phi)(x)|^2\,dx\right)\,d\nu(r_3),
\label{eq:covarsu2}
\ee
where we have replaced the $\th$-dependence by $r_3$-dependence, because
the map $\lm_{r_1,r_2}\,:G_{r_1,r_2}\to \Conj(SU(2))$ is injective, see
Eq.(\ref{eq:lmr1r2}).
\vspace*{.2cm}

If we now choose an explicit basis for the representation spaces 
we can explicitly calculate the Clebsch--Gordan
coefficients of ${\cal D}(SU(2))$. For the orthogonal bases we 
take the Wigner functions $D^j_{mn}$ as explained under Eq.(\ref{eq:Vrn}).

We will use the notation and definition of the Clebsch--Gordan coefficients
of $SU(2)$ as given in \cite{VMK}, chapter 8. Thus
\be
\D{j_1}{m_1 n_1}(g) \D{j_2}{m_2 n_2}(g) = \sum_{j=|j_1-j_2|}^{j_1+j_2} 
\sum_{m,n=-j}^j C^{jm}_{j_1 m_1 j_2 m_2} C^{jn}_{j_1 n_1 j_2 n_2}
\D{j}{mn}(g). 
\ee
The Clebsch--Gordan coefficients $C^{jm}_{j_1 m_1 j_2 m_2}$ are equal to 
zero if $m\neq m_1 + m_2$. So
\be
\D{j_1}{m_1 n_1}(g) \D{j_2}{m_2 n_2}(g) = \sum'_j 
C^{j\,(m_1+m_2)}_{j_1 m_1 j_2 m_2} C^{j\,(n_1+n_2)}_{j_1 n_1 j_2 n_2} 
\D{j}{(m_1+m_2)\,(n_1+n_2)}(g),
\ee
where the primed summation over $j$ runs from max($|j_1-j_2|, |m_1+m_2|,
|n_1+n_2|$) to $(j_1+j_2)$. 
In the tensor product representation $\Pi^{r_1}_{n_1}\otimes\Pi^{r_2}_{n_2}$
we consider the basis function 
\be
\Phi = D^{j_1}_{m_1 n_1}\otimes D^{j_2}_{m_2 n_2}\,:\,(y_1,y_2)\mapsto
 D^{j_1}_{m_1 n_1}(y_1)\,D^{j_2}_{m_2 n_2}(y_2),\quad j_i\geq n_i,\,
 -j_i\leq m_i\leq j_i, \,i=1,2
\ee
The mapping $\rho$ from Eq.(\ref{eq:rhon3}) takes this basis function to a 
linear combination of basis functions of a single irreducible  unitary 
representation $\Pi^{r_3}_{n_3}$:
\bea \label{eq:su2decom}
\left(\rho^{\th}_{n_3} \Phi\right)(x) &=& \int_{U(1)} e^{in_3\zeta}
\Phi(xg_{\zeta}w(\th)^{-1}, xg_{\zeta}w(\th)^{-1}a_{\th})\, d\zeta \nn\\
&=&\sum_{j=|j_1-j_2|}^{j_1+j_2}\sum_{m,p=-j}^{j} 
\sum_{p_2=-j_2}^{j_2}  C^{jm}_{j_1 m_1 j_2 m_2} C^{jp}_{j_1 n_1 j_2 p_2}
D^{j_2}_{p_2 n_2}(a_{\th})\,\times \nn\\
& &\qquad\qquad\int_{U(1)} e^{in_3\zeta} 
\sum_{r,s=-j}^{j} D^j_{mr}(x) D^j_{rs}(g_{\zeta}) D^j_{sp}(w(\th)^{-1})\,
d\zeta \\
&=& \sum'_j \left\{\sum'_{p_2} C^{j (m_1+m_2)}_{j_1 m_1 j_2 m_2} 
C^{j (n_1+p_2)}_{j_1 n_1 j_2 p_2} D^{j_2}_{p_2 n_2}
(a_{\th}) \overline{D^j_{(n_1+p_2) n_3}(w(\th))}\right\} \, D^j_{(m_1+m_2) 
n_3}(x)\nn
\eea
where the primed summation over $p_2$ runs from max($(-j-n_1), -j_2)$ to
min($(j-n_1), j_2)$. 

This shows how $\Phi \in V^{r_1}_{n_1}\otimes V^{r_2}_{n_2}$ can be
decomposed into single Wigner functions with a fixed label $n_3$, which
form a basis of $V^{r_3}_{n_3}$. 
The coefficients between the large brackets $\{\}$ now indeed are the
generalised Clebsch--Gordan coefficients for the quantum double group of
$SU(2)$. Clearly they depend on the representation labels, so on $(r_1,n_1),
(r_2,n_2)$ and $(r_3,n_3)$, where $r_3$ corresponds one--to--one to the
double coset $\th$. They also depend on the specific `states' labeled by 
the $j_1,m_1$, etc., just as one 
would expect. Note that $a_{\th}$ and $w(\th)$ 
are needed to implement the dependence on $\th$. We can denote 
these Clebsch--Gordan coefficients by
\be
\langle (r_1,n_1) j_1 m_1, (r_2,n_2) j_2 m_2\,|\, (r_3,n_3) j m\rangle := 
\sum'_{p_2} C^{jm}_{j_1 m_1 j_2 m_2} 
C^{j (n_1+p_2)}_{j_1 n_1 j_2 p_2} 
D^{j_2}_{p_2 n_2} (a_{\th}) \overline{D^j_{(n_1+p_2) n_3}(w(\th))}
\label{eq:clebgor}
\ee
with $r_3 = \lm_{r_1,r_2}(\th)$.
These coefficients are zero if $m\neq m_1+m_2$. Also, they are zero if
$n_3\neq (n_1+n_2)\;\mbox{mod}\;\dubbelZ$, so $n_3$ must be integer if 
$n_1 + n_2$ integer, and half integer if $n_1+n_2$ half integer. Thus we
can write
\be
\left(\rho^{r_3}_{n_3}\left(D^{j_1}_{m_1 n_1}\otimes D^{j_2}_{m_2 n_2}\right)
\right)(x) = \sum'_j \sum_{m=-j}^j \langle (r_1,n_1) j_1 m_1, (r_2,n_2) j_2 
m_2\,|\,(r_3,n_3) j m\rangle\,D^j_{m n_3}(x).
\ee

The isometry property of $\rho$ can now be calculated even more explicitly.
The left hand side of Eq.(\ref{eq:covarsu2}) gives 
\be
\int_{SU(2)}\int_{SU(2)} \D{j_1}{m_1 n_1}(y_1) \D{j_2}{m_2 n_2}(y_2)
\overline{\D{j_1}{m_1 n_1}(y_1)}\overline{\D{j_2}{m_2 n_2}(y_2)}\,dy_1 dy_2 =
\frac{1}{2j_1 + 1}\frac{1}{2j_2 +1}.
\label{eq:lhs}
\ee
For the right hand side of Eq.(\ref{eq:covarsu2}) we find
\be
\sum_{n_3}\int_{I_{r_1 r_2}}\left( \int_{SU(2)}|\rho^{\lm_{r_1 r_2}^{-1}
(r_3)}_{n_3} \left(\D{j_1}{m_1 n_1}
\otimes \D{j_2}{m_2 n_2}\right)(y)|^2\,dy\right)\,d\nu(r_3),
\ee
where $I_{r_1 r_2}$ given by Eq.(\ref{eq:imager1r2}), and the measure
$d\nu(r_3)$ by 
Eq.(\ref{eq:nur3}).  Substituting Eq.(\ref{eq:su2decom}) and 
Eq.(\ref{eq:clebgor}) yields
\bea
\lefteqn{\sum_{n_3}\int_{I_{r_1 r_2}} \int_{SU(2)}\left(\sum'_j\sum_m\langle 
(r_1,n_1) j_1 m_1, (r_2,n_2) j_2 m_2\,|\, (r_3,n_3) j m\rangle 
\D{j}{m n_3}(y)\right)\times\qquad}\nn\\
& &\left(\sum'_{j'}\sum_{m'} \overline{\langle (r_1,n_1) j_1 m_1, (r_2,n_2) 
j_2 m_2\,|\, (r_3,n_3) j' m'\rangle} \overline{\D{j'}{m'
n_3}(y)}\right) dy \,d\nu(r_3).
\eea
The integration over $y$ can be performed, and thus the isometry
property of the mapping $\rho$ reads
\be
\sum_{n_3}\int_{I_{r_1 r_2}} 
\sum'_j\frac{1}{2j+1}\,|\langle(r_1,n_1) j_1 m_1, (r_2,n_2) 
j_2 m_2\,|\,(r_3,n_3) j (m_1+m_2)\rangle|^2\,d\nu(r_3) = \frac{1}{2j_1
+1} \frac{1}{2j_2+1}.
\ee
More generally, if we start with the identity of inner products
which is immediately implied by Eq.(\ref{eq:covarsu2}), we obtain
\bea
\lefteqn{\sum_{n_3}\sum'_j \frac{1}{2j+1} \int_{I_{r_1 r_2}} \langle 
(r_1,n_1) j_1 m_1, (r_2,n_2) j_2 (m-m_1)\,|\, (r_3,n_3) j m\rangle 
\times\qquad}\nn\\
& &\overline{\langle (r_1,n_1) j'_1 m'_1, (r_2,n_2) 
j'_2 (m-m'_1)\,|\, (r_3,n_3) j m\rangle} \,d\nu(r_3) = \frac{\dl_{j_1,
j'_1}\,\dl_{j_2, j'_2}\,\dl_{m_1,m'_1}}{(2j_1+1)(2j_2+1)}
\label{eq:orthosu2}
\eea
This means that the Clebsch--Gordan coefficients (\ref{eq:clebgor})
for ${\cal D}(SU(2))$, built from Wigner functions and
Clebsch--Gordan coefficients for $SU(2)$, satisfy interesting
orthogonality relations, suggesting the existence of a `new' kind of
special functions. 

Remember that the $a_{\th}$ and $w(\th)$ given in Eqs.(\ref{eq:atheta}) and
(\ref{eq:nw}) are the choices we made for the 
Borel mappings $y(\xi)$ and $w(\xi)$ in Assumption \ref{ass2} and Proposition 
\ref{prop:maps} which uniquely depend on $r_3$ according to Eq.(\ref{eq:ys}). 
It is now clear that the choice of representatives in the double coset (so the
mapping $\xi \to y(\xi)$ of Assumption \ref{ass2}), and the choice of 
Borel map $\xi \to w(\xi)$ of Proposition \ref{prop:maps} 
do not affect the fusion rules: for $a_{\th} \mapsto g_{\phi}a_{\th}g_{\psi}$
and $w(\th)\mapsto g_{\phi} w(\th) g_{\zeta}$ the Clebsch--Gordan
coefficients from Eq.(\ref{eq:clebgor}) only change by a phase factor
$e^{i(n_1\phi -n_2\psi +n_3\zeta)}$, and thus the orthonormality
relations of Eq.(\ref{eq:orthosu2}) do not change.\\
This concludes our discussion of the fusion rules of ${\cal D}(SU(2))$. 
\section{Conclusion\label{s8}}
\setcounter{equation}{0}
In this paper we have focussed on the co-structure of the quantum
double $\DG$ of a compact group $G$ and have used it to study 
tensor products of irreducible representations. We have explicitly
constructed  a projection onto irreducible components for 
tensor product representations, which of course has to take into
account the (nontrivial) comultiplication. By subsequently using the 
Plancherel formula (i.e.\ by comparing squared norms) we found an
implicit formula for the multiplicities, or Clebsch--Gordan
series. Also, we have given the action of the universal $R$-matrix of
$\DG$ on tensor product states. For the example of $G=SU(2)$ we
calculated the Clebsch--Gordan series and coefficients
explicitly. In a forthcoming article we will expand further on the quantum
double of $SU(2)$, in particular the behaviour of its representations under
braiding and fusion. These results also will enable us to describe the 
quantum properties of topologically interacting point particles, as in
$\widetilde{ISO}(3)$ Chern--Simons theory, see \cite{BM}.
\appendix
\section{Some measure theoretical results\label{s9}}
\setcounter{equation}{0}
In this appendix we have collected some measure theoretical results
which have been used in section \ref{s6}.

\begin{thm}\label{111} (Kuratowski's theorem, see for instance
Parthasarathy, \cite{Par}, Ch.~I, Corollary 3.3)\\
If $E$ is a Borel subset of a complete separable metric space $X$
and $\lm$ is a one-one measurable map of $E$ into a separable metric space
$Y$ then $\lm(E)$ is a Borel subset of $Y$ and
$\lm\colon E\to\lm(E)$ is a Borel isomorphism.
\end{thm}

\begin{thm}\label{109} (Theorem of Federer \& Morse \cite{FedMor},
see also \cite{Par}, Ch.~I, Theorem 4.2)\\
Let $X$ and $Y$ be compact metric spaces and let $\lm$ be a continuous map of
$X$ onto $Y$. Then there is a Borel set $B\subset X$ such that $\lm(B)=Y$
and $\lm$ is one-to-one on $B$.
\end{thm}

The set $B$ is called a {\sl Borel section} for $\lm$.
Since the continuous image of a compact set is compact,
we can relax the conditions of Theorem \ref{109} by not requiring surjectivity
of $\lm$. Then $\lm(B)=\lm(X)$. By Theorem \ref{111} the mapping
$\lm|_B\colon B\to\lm(X)$ is a Borel isomorphism.
Let $\psi\colon\lm(X)\to B$ be the inverse of $\lm|_B$.
We will also call the mapping $\psi$ a {\sl Borel section} for $\lm$.
We conclude:

\begin{cor}\label{110}
Let $X$ and $Y$ be compact metric spaces and let $\lm$ be a continuous map of
$X$ to $Y$. Then there is a Borel map $\psi\colon\lm(X)\to X$ such that
$\lm(\psi(y))=y$ for all $y\in\lm(X)$ and $\psi(\lm(X))$ is a Borel set
in $X$.
\end{cor}

\begin{thm}\label{112} (isomorphism theorem, see for instance \cite{Par},
Ch.~I, Theorem 2.12)\\
Let $X_1$ and $X_2$ be two complete separable metric spaces and let
$E_1\subset X_1$ and $E_2\subset X_2$ be two Borel sets.
Then $E_1$ and $E_2$ are Borel isomorphic if and only if they have the
same cardinality.
In particular, if $E_1$ is uncountable, $X_2:=\dubbelR$ and $E_2$ is an open
interval, then $E_1$ and $E_2$ are Borel isomorphic.
\end{thm}

Next we discuss conditional probability, although we will not deal with
probabilistic interpretations. Our reference here is Halmos
\cite{Hal}, \S48.
Let $(X,\FSA)$ and $(Y,\FSB)$ be {\sl measurable spaces}, i.e.\ sets
$X$ and $Y$ with $\sg$-algebras $\FSA$ and $\FSB$, respectively.
Let $\lm\colon X\to Y$ be a measurable map.
Let $\mu$ be a probability measure on
$(X,\FSA)$. Define a probability measure $\nu$ on $(Y,\FSB)$ by the rule
\be
\nu(B):=\mu(\lm^{-1}(B)),\quad B\in\FSB.
\ee
By the Radon-Nikodym theorem there exists for each $A\in\FSA$
a $\nu$-integrable function $p_A$ on $Y$ such that
\be
\mu(A\cap \lm^{-1}(B))=\int_B p_A(y)\,d\nu(y),\quad B\in\FSB.
\ee
Then $p_A(y)$ is called the {\sl conditional probability of $A$ given $y$}.
Note that the functions $p_A$ are not unique. For fixed $A$, two choices for
$p_A$ can differ on a set of $\nu$-measure zero.
We will write
\be
p^y(A):=p_A(y),\quad y\in Y,\;A\in\FSA.
\ee
Then $p^y$ behaves in certain respects like a measure on $(X,\FSA)$,
but it may not be a measure.
If $f$ is a $\mu$-integrable function on $X$ then,
by the Radon-Nikodym theorem there exists a $\nu$-integrable function $e_f$
on $Y$ such that, for every $B\in\FSB$,
\be
\int_{\lm^{-1}(B)}f(x)\,d\mu(x)=\int_B e_f(y)\,d\nu(y).\label{114}
\ee

\begin{thm}\label{113}
If $\lm$ is a measurable map from a probability space
$(X,\FSA,\mu)$ to a measurable space $(Y,\nu)$, and if the conditional
probabilities $p_A(y)$ can be determined such that $p^y$ is a measure
on $(X,\FSA)$ for almost every $y\in Y$, then
\be
e_f(y)=\int_X f(x)\,dp^y(x)\quad\hbox{for $y$ almost everywhere on $Y$
w.r.t.\ $\nu$.}
\label{115}
\ee
In particular, if $X$ is an open interval in $\dubbelR$, or more generally
a complete separable metric space, then
$p_A(y)$ can be determined such that $p^y$ is a measure
on $(X,\FSA)$ for almost every $y\in Y$, and Eq.(\ref{114}) will hold
with $e_f(y)$ given by Eq.(\ref{115}).
\end{thm}

This theorem follows from Halmos \cite{Hal}, pp.\ 210--211, items (5) and (6)
together with the above Theorem \ref{112}.

Theorem \ref{113} greatly simplifies if $X$ is a complete separable metric
space and, moreover, $\lm$ is injective. Then
\bea
p^y(A) = p_A(y) = \chi_{\lm(A)}(y) = \left\{\begin{array}{cl} 0,& y\not\in
\lm(X),\nn\\ \dl_{\lm^{-1}(y)}(A),&y\in\lm(X) \end{array}\right.\nn\\
e_f (y) = \left\{\begin{array}{cl} 0, & y\not\in \lm(X),\nn\\f(\lm^{-1}(y)),
& y\in\lm(X),\end{array}\right.\nn\\
\int_{\lm^{-1}(B)} f(x)\,d\mu(x) = \int_{B\cap\lm(X)} f(\lm^{-1}(y))\,d\nu(y)
\eea 

\subsection*{Acknowledgements}
The third author was supported by the Dutch Science Foundation FOM/NWO.
We would like to thank Dr A.A.\ Balkema for useful discussions on
measure theory. Also, we want to thank dr. E.P. van den Ban for his
private communication 
\cite{vdBan}, which we used in discussing our Assumption \ref{asssplit}.

\end{document}